\documentclass[11pt]{article}

\usepackage[top=2.54cm, bottom=2.54cm, left=2.5cm, right=2.5cm]{geometry}
\usepackage{graphicx}
\graphicspath{ {./images/} }
\usepackage[backend=bibtex,style=alphabetic,backref]{biblatex}
\usepackage{titletoc}
\usepackage[big]{titlesec}
\usepackage{amsmath,amssymb,amsfonts,amsthm}
\usepackage{MnSymbol}
\usepackage{mathtools}
\usepackage{enumitem}
\usepackage{prftree}
\usepackage{turnstile}
\usepackage[cal=euler]{mathalpha}

\usepackage{hyperref}
\hypersetup{colorlinks=true,allcolors=DarkBlue}
\newcommand*{\fullref}[1]{\hyperref[{#1}]{\nameref*{#1} \ref*{#1}}}
\newcommand*{\secref}[1]{\hyperref[{#1}]{\S\ref*{#1} \nameref*{#1}}}

\usepackage{tikz}
\usetikzlibrary{cd}

\title{Kreisel--L{\'e}vy-type theorems for Kripke--Platek and other set theories}

\author{Shuangshuang Shu, Michael Rathjen}

\date{}

\titleformat{name=\subsection}[block]
   {\normalfont\bfseries}{\thesubsection}{1em}{}
\titleformat{name=\subsection,numberless}
   {\vspace{.8ex}\normalfont\scshape\large}{}{0pt}{}

\definecolor{Red}{rgb}{0.70,0.00,0.10}
\definecolor{Cyan}{rgb}{0.00,0.20,0.55}
\definecolor{DarkGreen}{rgb}{0.00,0.20,0.00}
\definecolor{Blue}{rgb}{0.00,0.10,0.70}
\definecolor{DarkBlue}{rgb}{0.00,0.05,0.30}

\setlist[enumerate,1]{label={(\roman*)},nosep}
\setlist[enumerate,2]{nosep}
\setlist[itemize,1]{nosep}
\setlist[itemize,2]{nosep}

\numberwithin{equation}{section}

\newtheoremstyle{defaultstyle}
{\topsep}{\topsep}
{\normalfont}{}
{\bfseries}{.}
{0.5em}
{\thmname{#1}\thmnumber{ #2}}

\newtheoremstyle{customstyle}
{\topsep}{\topsep}
{\normalfont}{}
{\bfseries}{.}
{0.5em}
{\thmname{#3}\thmnumber{ #2}}

\newtheoremstyle{customremarkstyle}
{\topsep}{\topsep}
{\normalfont}{}
{\itshape}{.}
{0.5em}
{\thmname{#3}\thmnumber{ #2}}

\theoremstyle{defaultstyle}
\newtheorem{thm1}{Theorem}[section]
\newtheorem{def1}[thm1]{Definition}
\newtheorem{prop1}[thm1]{Proposition}
\newtheorem{cor1}[thm1]{Corollary}
\newtheorem{lem1}[thm1]{Lemma}
\newtheorem{ex1}[thm1]{Exercise}

\theoremstyle{customstyle}
\newtheorem{named}[thm1]{NAME}
\newtheorem*{named*}{NAME}

\theoremstyle{remark}
\newtheorem{obsdep}{Observation}[thm1]

\theoremstyle{customremarkstyle}
\newtheorem{namedrmk}[thm1]{NAME}
\newtheorem*{namedrmk*}{NAME}

\newtheorem{sublemma}{NAME}[thm1]

\newenvironment{thm}{\begin{thm1}[Theorem]}{\end{thm1}}
\newenvironment{defn}{\begin{def1}[Definition]}{\end{def1}}
\newenvironment{prop}{\begin{prop1}[Proposition]}{\end{prop1}}
\newenvironment{cor}{\begin{cor1}[Corollary]}{\end{cor1}}
\newenvironment{lem}{\begin{lem1}[Lemma]}{\end{lem1}}

\makeatletter
\def\blfootnote{\xdef\@thefnmark{}\@footnotetext}
\makeatother

\let\leq\leqslant
\let\geq\geqslant
\let\emptyset\varnothing
\let\epsilon\varepsilon

\mathchardef\mhyphen="2D

\newcommand{\PA}{\mathrm{PA}}

\newcommand{\True}{\mathrm{True}}

\newcommand{\KP}{\mathrm{KP}}
\newcommand{\ON}{\mathrm{ON}}
\newcommand{\TI}{\mathrm{TI}}

\newcommand{\Prog}{\mathrm{Prog}}
\newcommand{\Prov}{\mathrm{Prov}}

\newcommand{\bV}{\mathbf{V}}

\newcommand{\RFN}{\mathrm{RFN}}

\newcommand{\RS}{\mathrm{RS}}

\newcommand{\RSderivation}{\mathrm{RS}_\Omega(\bV)\mhyphen\mathrm{derivation}}
\newcommand{\RSquasicode}{\mathrm{RS}_\Omega(\bV)\mhyphen\mathrm{quasicode}}
\newcommand{\FinSeq}{\mathrm{Set}^{<\omega}}

\newcommand{\convg}[1]{{{#1}\!\!\downarrow}}

\newcommand{\qedbyhand}{\hspace*{\fill}\qedsymbol}

\DeclareMathOperator{\dom}{dom}

\DeclareMathOperator{\rank}{rank}

\DeclarePairedDelimiter\code{\ulcorner}{\urcorner}

\begin{document}


\maketitle

\blfootnote{Department of Pure Mathematics, University of Leeds, Leeds, LS2 9JT, United Kingdom\\
\hspace*{1.8em}mmssh@leeds.ac.uk, m.rathjen@leeds.ac.uk}

\begin{abstract}
We prove that, over Kripke--Platek set theory with infinity (KP), transfinite induction along the ordinal $\epsilon_{\Omega+1}$ is equivalent to the schema asserting the soundness of KP, where $\Omega$ denotes the supremum of all ordinals in the universe; this is analogous to the result that, over Peano arithmetic (PA), transfinite induction along $\epsilon_0$ is equivalent to the schema asserting the soundness of PA. In the proof we need to code infinitary proofs within KP, and it is done by using partial recursive set functions. This result can be generalised to $\KP+\Gamma\text{-separation}+\Gamma\text{-collection}$ where $\Gamma$ is any given syntactic complexity, but not to ZF. 
\end{abstract}

\setlength{\parindent}{10pt}

\tableofcontents

\pagebreak

\section{Introduction}

In 1968, \cite{KreiselLevy1968} showed that over Peano arithmetic (PA), the axiom schema which asserts the soundness of PA is equivalent to the schema of transfinite induction along $\epsilon_0$, the first ordinal $\alpha$ that satisfies $\alpha=\omega^\alpha$: let 
\begin{align*}
\RFN(\PA)=\{\forall x(\mathrm{Prov}_{\PA}(\code{\phi(\dot{x})})\to \phi(x))\mid \phi\text{ a formula in PA}\}
\end{align*}
be the soundness principle of PA\footnote{``RFN'' represents ``reflection'' which is the name of the schema used in \cite{KreiselLevy1968}; in set theory, the same phrase can mean something else, so we avoid using this term}, where $\mathrm{Prov}_{\PA}(\code{\phi(\dot{x})})$ means $\phi(x)$ is $\PA$-provable\footnote{we assume a standard coding $\code{\cdot}$ for formulas; $\code{\phi(\dot{x})}$ is a primitive-recursive function which on input $m$ outputs $\code{\phi(S^m(0))}$ where $S$ is the symbol for successor}; and let 
\begin{align*}
\epsilon_0=\sup_n\omega_n\quad\text{where }\omega_0=\omega,\,\omega_{n+1}=\omega^{\omega_n};
\end{align*}
the transfinite induction schema along an ordinal $\alpha$ is 
\begin{align*}
\TI(\alpha)=\{ \forall \beta((\forall \gamma<\beta)\phi(\gamma)\to \phi(\beta))\to (\forall \beta<\alpha)\phi(\beta) \mid \phi\text{ a formula in PA} \};
\end{align*}
then the theorem mentioned above reads as 
\begin{align*}
\RFN(\PA)\equiv \TI(\epsilon_0)\quad\text{over PA},
\end{align*}
i.e., $\PA+\RFN(\PA)$ proves $\TI(\epsilon_0)$ and vice versa. The direction $\PA+\RFN(\PA)\vdash \TI(\epsilon_0)$ is essentially from Gentzen's proof that $\PA$ proves $\TI(\alpha)$ for every $\alpha<\epsilon_0$, and \cite{KreiselLevy1968} used another device to prove the other direction; in \cite{Schwichtenberg1977}, the direction $\PA+\TI(\epsilon_0)\vdash \RFN(\PA)$ is proved by transfinite induction along infinitary proof trees as which all PA-proofs can be interpreted with lengths $<\epsilon_0$. \par 

This paper aims to prove a very similar result in Kripke--Platek set theory (KP\footnote{contrary to convention, we will use KP to denote Kripke--Platek set theory with infinity}). KP is a subtheory of ZF consisting of the following axioms:
\begin{enumerate}[label=\arabic*.]
\item Extensionality. $\forall x(x\in a\leftrightarrow x\in b)\to a=b$.
\item Set induction. $\forall x((\forall y\in x)G(y)\to G(x))\to \forall x\,G(x)$, for all formulas $G$.
\item Pair. $\exists x(x=\{ a,b\})$.
\item Union. $\exists x(x=\bigcup a)$.
\item Infinity. $\exists x(x\neq \emptyset \wedge (\forall y\in x)(\exists z\in x)(y\in z))$.
\item $\Delta_0$-separation. $\exists x\forall u(u\in x\leftrightarrow (u\in a\wedge F(u)))$, for all $\Delta_0$-formulas $F$.
\item $\Delta_0$-collection. $(\forall x\in a)\exists y\,G(x,y)\to \exists z(\forall x\in a)(\exists y\in z)G(x,y)$, for all $\Delta_0$-formulas $G$. 
\end{enumerate}
In proof theory, KP is rather similar to PA. For example, it is a well-known result in \cite{Parsons1970} that the primitive-recursive functions on natural numbers are exactly the class of the provably total functions of PA with induction restricted to $\Sigma_1$-formulas. In set theory we may define primitive-recursive functions on sets as well, and by \cite{Rathjen1992} Theorem 1.2, a set function $F$ is primitive-recursive iff it is provably total in KP with (set) induction restricted to $\Sigma_1$-formulas. But more importantly for us, there is a way to transform KP-proofs into infinitary proof trees for which the cuts can at least be partially eliminated. Thus, if we let $\Omega$ denote the supremum of all ordinals in the universe of KP, and propose that 
\begin{align*}
\RFN(\KP)\equiv \TI(\epsilon_{\Omega+1})\quad\text{over KP}
\end{align*}
where $\epsilon_{\Omega+1}$ denotes the $(\Omega+1)$-th ordinal $\alpha$ that satisfies $\alpha=\omega^\alpha$, then it is possible to follow the method mentioned in \cite{Schwichtenberg1977} to prove our proposition; the equivalence stated above is what we are trying to prove in this paper. \par

In Section \ref{sec:Preliminary definitions} we lay out some preliminary definitions. The direction $\KP+\RFN(\KP)\vdash \TI(\epsilon_{\Omega+1})$ will be quickly dealt with in \ref{sec:Provable well-orderings in KP}, and all the rest \ref{sec:The infinitary proof system}--\ref{sec:Embedding KP} are dedicated to proving the converse $\KP+\TI(\epsilon_{\Omega+1})\vdash \RFN(\KP)$. The length of the latter is partly due to the lengthiness of cut elimination theorem and embedding theorem themselves, but also due to a subtlety in expressing the infinitary proofs within KP, and the need of using partial recursive set functions for which the recursion theorem applies.

\section{Preliminary definitions}
\label{sec:Preliminary definitions}

In this paper, our theory $T$ is always assumed to be sufficiently strong, say, it contains primitive-recursive arithmetic.

\begin{defn}
Let $T$ be a theory. The \textbf{soundness principle} of $T$, $\RFN(T)$, is the schema
\begin{align*}
\RFN(T)=\{\forall x(\mathrm{Prov}_{T}(\code{\phi(\dot{x})})\to \phi(x))\mid \phi\text{ a formula in }T\}
\end{align*}
where $\mathrm{Prov}_{T}(\code{\phi(\dot{x})})$ means that $\phi(x)$ is $T$-provable.
\end{defn}

\begin{defn}
\label{primitive-recursive set function}
If our background theory is a set theory in which $\omega$ is definable, \textbf{primitive-recursive set functions}\footnote{the initial function $x\mapsto \omega$ is dispensable in a more general setting} are the class of functions built up from the initial functions
\begin{itemize}
\item $P_{n,i}(x_1,\ldots,x_n)=x_i$, $1\leq i\leq n$, 
\item $Z(y)=0$,
\item $M(x,y)=x\cup\{ y\}$,
\item $C(x,y,u,v)=x$ if $u\in v$, $y$ otherwise,
\item $\mathrm{Inf}(x)=\omega$
\end{itemize}
by 
\begin{itemize}
\item \emph{substitution}: $F(\vec{x})=K(G_1(\vec{x}),\ldots,G_k(\vec{x}))$, and
\item \emph{primitive recursion}: $F(z,\vec{x})=H(\bigcup\{ F(u,\vec{x})\mid u\in z\},z,\vec{x} )$. 
\end{itemize}
\end{defn}

Some of the most important examples of primitive-recursive set functions are the characteristic functions of $\Delta$-predicates in KP, and ordinal operations such as addition, multiplication and exponentiation. See \cite{Rathjen1992} 2.2 for details.

\begin{defn}
Let $T$ be a theory. An \textbf{ordinal representation system} $\langle R,\preceq\rangle$ in $T$ is a unary relation $R$ and a preorder $\preceq$ primitive-recursively coded in $T$; members of $R$ are strings of symbols to be interpreted as ordinals and $\preceq$ orders the ordinals; the system also comes with basic ordinal operations: addition, multiplication, exponentiation $\alpha\mapsto\omega^\alpha$, all primitive-recursively coded in $T$. 
\end{defn}

The strings of symbols in $R$ may be mapped non-injectively into the ordinals and so $\preceq$ could be not antisymmetric. For example, $\omega$ and $\omega^1$ have the same value, $\omega\preceq \omega^1$ and $\omega^1\preceq \omega$, but their notations are different: $\omega\not\equiv \omega^1$.

\begin{defn}
Let $T$ be a theory and let $\langle R,\preceq\rangle$ be an ordinal representation system in $T$; let $U$ be a predicate. The \textbf{progressiveness of $\prec$ for $U$}, $\Prog_\prec(U)$, is the formula
\begin{align*}
\Prog_\prec(U)\equiv \forall x((\forall y\prec x)y\in U\to x\in U ).
\end{align*}
If $a\in R$, the \textbf{transfinite induction along $a$ for $U$}, $\TI_\prec(a,U)$, is the formula
\begin{align*}
\TI_\prec(a,U)\equiv \Prog_\prec(U)\to (\forall x\prec a)x\in U.
\end{align*}
If $F$ is a formula, then by $x\in F$ we mean $F(x)$; so, e.g. $\Prog_\prec(F)$ is the formula $\forall x((\forall y\prec x)F(y)\to F(x) )$. The \textbf{transfinite induction schema along $a$}, $\TI_\prec(a)$, is the schema
\begin{align*}
\TI_\prec(a)=\{ \TI_\prec(a,F)\mid F\text{ a formula in }T\}.
\end{align*}
\end{defn}

If it is clear which ordinal representation system we are using, we often drop the $\prec$ symbol. \par

The ordinal representation system in KP we are going to use is rather canonical. In arithmetic, apparently the numbers we have access to are the natural numbers, but by Cantor normal form, every non-zero ordinal $\alpha$ can be written in the form
\begin{align*}
\alpha=\omega^{\alpha_1}+\ldots+\omega^{\alpha_n}
\end{align*} 
for some $n\in \omega$ and ordinals $\alpha_1\geq\ldots\geq \alpha_n$; if $\alpha<\epsilon_0$ then $\alpha_1,\ldots,\alpha_n$ must be smaller than $\alpha$, and we can write $\alpha_1,\ldots,\alpha_n$ in Cantor normal forms again, representing them by smaller ordinals. The descending sequences of ordinals must end in finitely many steps, therefore we may represent every ordinal $\alpha<\epsilon_0$ by just natural numbers. The idea for representing $\epsilon_{\Omega+1}$ in KP is exactly the same; we will only sketch the definition.

\begin{defn}
Within KP we define an ordinal representation system $\langle R,\preceq\rangle$ as follows. Let ON denote the class of ordinal numbers. If $a,b$ are members of $R$, we use $a\equiv b$ to denote that $a,b$ are exactly the same string of symbols, whereas $a=b$ means the weaker condition $a\preceq b$ and $b\preceq a$.
\begin{itemize}
\item The symbols are $\{\alpha\mid \alpha\in \ON\}$, $\Omega$, $+$, $\omega^\cdot$, $\epsilon_\cdot$.
\item The members $a$ of the class $R$ are in the following forms. 
\begin{itemize}
\item $a\equiv \alpha$ for some $\alpha\in \ON$. 
\item $a\equiv \Omega$. 
\item $a\equiv \omega^{a_0}$ where $a_0\in R$. 
\item $a\equiv a_1+a_2$ where $a_1,a_2$ are of the form $\omega^{b_1}+\ldots+\omega^{b_n}$ (not necessarily in Cantor normal form) where $b_1,\ldots,b_n\in R$. 
\item $a\equiv \epsilon_{a_0}$ where $a_0\in R$. 
\end{itemize}
\item The ordering $\prec$ is defined as follows. Let $a,b\in R$. 
\begin{itemize}
\item $a,b\in \ON$. Then $a\prec b$ iff $a<b$ as ordinals. 
\item $a\in \ON$, $b\equiv\Omega$. Then $a\prec b$. 
\item $a\in \ON$, $b\equiv b_1+b_2$ or $b\equiv \omega^{b_0}$ or $b\equiv \epsilon_{b_0}$. \par 
If $b$ contains any $\Omega$ symbol, then $a\prec b$. \par 
If $b$ does not contain any $\Omega$ symbol, then $b$ can be evaluated in ON primitive-recursively, and $a\prec b$ iff $a<b$ in the evaluation. 
\item $a\equiv\Omega$, $b\equiv b_1+b_2$. \\
If $b_1\prec \Omega$, then $b\diamond \Omega$ iff $b_2\diamond\Omega$, where $\diamond$ is any of $\prec,=,\succ$. \\
If $b_1=\Omega$, then $b\succ \Omega$ (since $b_2\succ 0$ always in our convention for $+$). \\
If $b_1\succ\Omega$, then $b\succ \Omega$. 
\item $a\equiv\Omega$, $b\equiv \omega^{b_0}$ or $\epsilon_{b_0}$. Then $b\diamond \Omega$ iff $b_0\diamond \Omega$ for $\diamond\in \{ \prec,=,\succ\}$. 
\item $a\equiv a_1+a_2$, $b\equiv b_1+b_2$ or $b\equiv \omega^{b_0}$. We may primitive-recursively sort $a$ and $b$ into their Cantor normal form and then compare them in the usual way. 
\item $a\equiv a_1+a_2$, $b\equiv \epsilon_{b_0}$. Let $\omega^{a_0}$ be the first term of the Cantor normal form of $a$. \\
If $a_0\prec b$, then $a\prec b$. \\
If $a_0\succeq b$, then $a\succ b$.
\item $a\equiv \omega^{a_0}$, $b\equiv \omega^{b_0}$. Then $a\diamond b$ iff $a_0\diamond b_0$ for $\diamond\in \{ \prec,=,\succ\}$. 
\item $a\equiv \omega^{a_0}$, $b\equiv \epsilon_{b_0}$. Then $a\diamond b$ iff $a_0\diamond b$ for $\diamond\in \{ \prec,=,\succ\}$.
\end{itemize}
\end{itemize}
The predicates $R$ and $\prec$ are therefore $\Delta$ in $\KP$ and are primitive-recursive.
\end{defn}


We will be using an auxiliary sum operation for ordinals.

\begin{defn}
\label{natural sum definition}
Define an operation $\#$ as follows: $\alpha\#0=0\#\alpha=\alpha$. For $\alpha$ with Cantor normal form $\omega^{\alpha_1}+\ldots+\omega^{\alpha_m}$ and $\beta$ with Cantor normal form $\omega^{\alpha_{m+1}}+\ldots+\omega^{\alpha_n}$ define
\begin{align*}
\alpha\#\beta:= \omega^{\alpha_{\pi(1)}}+\ldots +\omega^{\alpha_{\pi(n)}}
\end{align*}
where $\pi$ is a permutation of $\{ 1,\ldots,n\}$ such that $\alpha_{\pi(i)}\geq \alpha_{\pi(i+1)}$ for all $1\leq i<n$. $\alpha\#\beta$ is called the \textbf{natural sum} or \textbf{Hessenberg sum} of $\alpha$ and $\beta$.
\end{defn}

\section{Provable well-orderings in \texorpdfstring{$\KP$}{KP}}
\label{sec:Provable well-orderings in KP}

In this section we prove the direction $\KP+\RFN(\KP)\vdash \TI(\epsilon_{\Omega+1})$; the proof follows almost directly from modification of Gentzen's proof that PA proves $\TI(\alpha)$ for every $\alpha<\epsilon_0$. \par

Recall that $R$ denotes the symbol class for our ordinal representation system and $\prec$ is our ordering on $R$.

\begin{defn}
For a predicate $U$, let 
\begin{align*}
U^J=\{ b\in R\mid (\forall a\in R)(R\cap a\subset U\to R\cap a+ \omega^b\subset U)\},
\end{align*}
denote the \textbf{jump} of $U$, where $R\cap x$ means $\{ b\in R\mid b\prec x\}$. 
\end{defn}

\begin{lem}
$\KP\vdash \Prog_\prec(U)\to \Prog_\prec(U^J)$. 
\end{lem}
\begin{proof}
Assume (1) $\Prog_\prec(U)$ and (2) $(\forall x\prec b)x\in U^J$, we want to show that $b\in U^J$, i.e., 
\begin{align*}
(\forall a\in R)(R\cap a\subset U\to R\cap a+ \omega^b\subset U).
\end{align*}
Assume that (3) $R\cap a\subset U$. Let $d\in R\cap a+\omega^b$; we need to show that $d\in U$ under the assumptions (1)--(3). \par 
If $d\prec a$, then $d\in U$ by (3). \par 
If $d=a$ then using (1) and (3) we have $d\in U$. \par 
If $d\succ a$, then since $d\prec a+\omega^b$, we may primitive-recursively find $d_1,\ldots,d_k$ such that 
\begin{align*}
d=a+\omega^{d_1}+\ldots+\omega^{d_k},\quad d_k\preceq\ldots\preceq d_1\prec b.
\end{align*}
Since $R\cap a\subset U$, we get $R\cap a+\omega^{d_1}\subset U$ by (2). Using (2) a further $k-1$ times we obtain 
\begin{align*}
R\cap a+\omega^{d_1}+\ldots+\omega^{d_k}\subset U.
\end{align*}
Finally, using one application of (1) we have $d\in U$. 
\end{proof}

\begin{namedrmk*}[Notation]
Let $e_0\equiv \Omega+1$, $e_{n+1}\equiv \omega^{e_n}$.
\end{namedrmk*}

The following lemma shows that $\KP\vdash \TI_\prec(e_n,U)$ for all $n\in \omega$.

\begin{lem}
\label{Gentzen argument in KP}
For any $n<\omega$ and any definable class $U$, 
\begin{align*}
\KP\vdash \Prog_\prec(U)\to R\cap e_n\subset U\wedge e_n\in U.
\end{align*}
\end{lem}
\begin{proof}
By induction on $n$ (outside of KP). \par 
For $n=0$ we need to show that $\forall x((\forall y\prec x)y\in U\to x\in U)\to \Omega+1\subset U\wedge \Omega+1\in U$. Suppose $\Prog_\prec(U)$ holds but there is some $a\prec \Omega$ such that $a\notin U$. Then by $\in$-induction of KP, there is a least ordinal $\alpha$ such that there exists some $b\in R$ with value $\alpha$ and $b\notin U$. But this implies $(\forall y\prec b)y\in U$, and therefore $b\in U$ by $\Prog_\prec(U)$. Thus, assuming $\Prog_\prec(U)$, we must have $\Omega\subset U$. Then by applying $\Prog_\prec(U)$ twice we have $\Omega,\Omega+1\in U$. \par 
Now suppose the result holds for $n$. Since the induction hypothesis (i.h.) holds for all definable classes, we have that 
\begin{align*}
\KP\vdash \Prog_\prec(U^J)\to R\cap e_n\subset U^J\wedge e_n\in U^J.
\end{align*}
Since $\KP\vdash \Prog_\prec(U)\to \Prog_\prec(U^J)$, we have 
\begin{align*}
(*)\quad \KP\vdash \Prog_\prec(U)\to R\cap e_n\subset U^J\wedge e_n\in U^J. 
\end{align*}
Now we argue in KP. Assume $\Prog_\prec(U)$, then from $(*)$, we obtain 
\begin{align*}
R\cap e_n\subset U^J\wedge e_n\in U^J. 
\end{align*}
By the definition of $U^J$, $e_n\in U^J$ implies that $R\cap 0\subset U \to R\cap 0+\omega^{e_n}\subset U$. Thus $R\cap \omega^{e_n}\subset U$, and an application of $\Prog_\prec(U)$ yields $\omega^{e_n}\in U$ as required. 
\end{proof}

\begin{lem}
$\KP\vdash \forall n(\Prov_{\KP}(\code{\TI_\prec(e_{\dot{n}},F)}))$ for any formula $F$.
\end{lem}
\begin{proof}
Let $F$ be given; we describe the procedure of writing the proof for $\TI_\prec(e_n,F)$. For any definable class $U$, let $U^{J^n}$ denote $U^{J\ldots J}$ with $J$ applied $n$ times. To find the proof for $\TI_\prec(e_n,F)$, we start with the proof of $\Prog_{\prec}(F^{J^n})\to R\cap e_0\subset F^{J^n}\wedge e_0\in F^{J^n}$. Applying the argument of the previous lemma \ref{Gentzen argument in KP}, we obtain a proof of $\Prog_\prec(F^{J^{n-1}})\to R\cap e_1\subset F^{J^{n-1}}\wedge e_1\in F^{J^{n-1}}$. Then the argument can be applied again to obtain a proof of $\Prog_\prec(F^{J^{n-2}})\to R\cap e_2\subset F^{J^{n-2}}\wedge e_2\in F^{J^{n-2}}$, and so on. Eventually, we arrive at $\Prog_\prec(F)\to R\cap e_n\subset F\wedge e_n\in F$, implying $\TI_\prec(e_n,F)$. This proof-writing function with argument in $n$ is primitive-recursive, hence provably total in KP, which implies that $\KP\vdash \forall n(\Prov_{\KP}(\code{\TI_\prec(e_{\dot{n}},F)}))$.
\end{proof}


\begin{cor}
$\KP+\RFN(\KP)\vdash \forall n\,\TI_\prec(e_n,F)$.  \qedbyhand 
\end{cor}

\begin{lem}
$\KP\vdash \forall x(x\prec \epsilon_{\Omega+1}\leftrightarrow (\exists n\in \omega)x\prec e_n )$. 
\end{lem}
\begin{proof}
Clearly $\KP\vdash \forall x((\exists n\in \omega)x\prec e_n\to x\prec \epsilon_{\Omega+1})$. The other direction is immediate for the cases $a\in \ON$, $a\equiv \Omega$, or $a\equiv \epsilon_{a_0}$. The cases $a\equiv a_1+a_2$ or $a\equiv \omega^{a_0}$ are done by induction on the lengths of expressions in $R$. Let $\omega^{a_0}$ be the first term of the Cantor normal form of $a$. If $a\prec \epsilon_{\Omega+1}$, we must have $a_0\prec \epsilon_{\Omega+1}$. Then by the i.h., $a_0\prec e_n$ for some $n$ and hence $\omega^{a_0}\prec e_{n+1}$, so $a\prec e_{n+1}$.
\end{proof}

From this lemma we obtain

\begin{cor}
\label{easy direction of main theorem}
$\KP+\RFN(\KP)\vdash \TI_\prec(\epsilon_{\Omega+1},F)$ for any formula $F$.   \qedbyhand 
\end{cor}

\section{The infinitary proof system}
\label{sec:The infinitary proof system}

The rest of the paper is dedicated to prove 
\begin{align*}
\KP+\TI(\epsilon_{\Omega+1})\vdash \RFN(\KP). 
\end{align*}
This will be done by considering KP-proofs as infinitary proofs with lengths $<\epsilon_{\Omega+1}$ and with cut complexity at most $\Pi_1/\Sigma_1$; an induction along such proof trees will show that the proofs have true conclusions, thus fulfilling $\RFN(\KP)$. \par

We first introduce our infinitary proof system $\RS_\Omega(\bV)$, which is a Tait-style sequent calculus.

\begin{named}[Definitions of $\RS_\Omega(\bV)$-terms and -formulas]
\label{definition of RS(V)-terms and -formulas}
$ $\par 
\begin{itemize}
\item For every set $a$, the constant $c_a$ is an $\RS_\Omega(\bV)$-term.
\item If $s,t$ are $\RS_\Omega(\bV)$-terms, so are $\{ s,t\}$ and $\bigcup s$. 
\item If $s,t_1,\ldots,t_n$ are $\RS_\Omega(\bV)$-terms and $A(a,b_1,\ldots,b_n)$ is a $\Delta_0$-formula of KP with all free variables displayed, then $\{ x\in s\mid A(x,t_1,\ldots,t_n)\}$ is an $\RS_\Omega(\bV)$-term.
\item If $t_1,\ldots,t_n$ are $\RS_\Omega(\bV)$-terms and $A(a_1,\ldots,a_n)$ is a formula of KP with all free variables displayed, then $A(t_1,\ldots,t_n)$ is an $\RS_\Omega(\bV)$-formula.
\end{itemize}
\end{named}

$\RS_\Omega(\bV)$-terms are purely symbolic, but we can still evaluate them as sets. Let Comp be a primitive-recursive set function that does the following: for all sets $s,t_1,\ldots,t_n$, $\Delta_0$-formula $A(a,b_1,\ldots,b_n)$ with all free variables indicated, 
\begin{align*}
\mathrm{Comp}( s,\vec{t},\code{A})=\{ x\in s\mid A(x,\vec{t})\};
\end{align*}  
a primitive-recursive computation of such a function follows from, e.g., the proof of \cite{Barwise1975} I.5.2 (v). Then we can define a primitive-recursive set function Ev which evaluate $\RS_\Omega(\bV)$-terms in the set universe:
\begin{itemize}
\item $\mathrm{Ev}(u)=a$ if $u\equiv c_a$ is a constant for the set $a$;
\item if $u\equiv \bigcup s$, then $\mathrm{Ev}(u)=\bigcup\mathrm{Ev}(s)$;
\item if $u\equiv \{ s,t\}$, then $\mathrm{Ev}(u)=\{ \mathrm{Ev}(s),\mathrm{Ev}(t)\}$;
\item if $u\equiv \{ x\in s\mid A(x,t_1,\ldots,t_n)\}$, then $\mathrm{Ev}(u)=\mathrm{Comp}( s,\vec{t},\code{A})$. 
\end{itemize}



\begin{namedrmk*}[Notation]
$ $\par 
\begin{itemize}
\item If $u$ is an $\RS_\Omega(\bV)$-term, $|u|$ denotes the set-theoretic rank of $\mathrm{Ev}(u)$. 
\item The formula $s=t$ is a shorthand for $(\forall x\in s)x\in t\wedge (\forall x\in t)x\in s$.
\item If $A$ is a formula, $\pm A$ denotes formulas both $A$ and $\neg A$. 
\item If $A$ is a formula, and $z$ is a variable not appearing in $A$, then $A^z$ denotes $A$ relativised to $z$: $A^z$ is the result of replacing every unbounded quantifiers $\exists x$ by $(\exists x\in z)$ and $\forall x$ by $(\forall x\in z)$ in $A$. 
\end{itemize}
\end{namedrmk*}

\begin{named}[Definition of derivability in $\RS_\Omega(\bV)$]
We give an inductive definition of the relation $\RS_\Omega(\bV) \sststile{}{\alpha}\Gamma$ by recursion on $\alpha$; this definition is yet outside of KP. The symbol $\Gamma$ stands for an arbitrary finite set of $\RS_\Omega(\bV)$-formulas, and if $A$ is a formula, $\Gamma,A$ means $\Gamma\cup\{ A\}$. $\RS_\Omega(\bV)\sststile{}{\alpha} \Gamma$ is meant to express that the system $\RS_\Omega(\bV)$ proves the disjunction $\bigvee\Gamma$ with a proof of length $\preceq\alpha$. \par 
The axioms of $\RS_\Omega(\bV)$ are of the form $\sststile{}{\alpha}\Gamma,A$ where $A$ is a $\Delta_0$-formula true in KP (i.e. provable in KP). More precisely, if $A\equiv A(t_1,\ldots,t_n)$ is $\Delta_0$, $t_1,\ldots,t_n$ are all the $\RS_\Omega(\bV)$-terms appearing in the formula, and $A(\mathrm{Ev}(t_1),\ldots,\mathrm{Ev}(t_n))$ is true in KP, then $\Gamma,A$ is an axiom of $\RS_\Omega(\bV)$. \par 
The following are the inference rules of $\RS_\Omega(\bV)$. 
\begin{align*}
&(\wedge)\quad \prftree{ \sststile{}{\alpha_0}\Gamma,A\quad \sststile{}{\alpha_1}\Gamma,B}{\sststile{}{\alpha}\Gamma,A\wedge B}\quad \alpha_0,\alpha_1<\alpha\\
&(\vee) \quad \prftree{\sststile{}{\alpha_0}\Gamma,C\quad\text{for some }C\in \{ A,B\}}{\sststile{}{\alpha}\Gamma,A\vee B}\quad \alpha_0<\alpha\\
& (b\forall) \quad \prftree{\sststile{}{\alpha_s}\Gamma,s\in t\to A(s)\quad\text{for all terms }s}{\sststile{}{\alpha}\Gamma,(\forall x\in t)A(x)}\quad \alpha_s<\alpha\\
& (b\exists) \quad \prftree{\sststile{}{\alpha_0}\Gamma,s\in t\wedge A(s)}{\sststile{}{\alpha}\Gamma,(\exists x\in t)A(x)}\quad \alpha_0<\alpha \\  
\end{align*}
\begin{align*}
& (\forall)\quad \prftree{\sststile{}{\alpha_s}\Gamma,A(s)\quad\text{for all terms }s}{\sststile{}{\alpha}\Gamma,\forall x\,A(x)}\quad \alpha_s<\alpha \\
& (\exists )\quad \prftree{\sststile{}{\alpha_0}\Gamma,A(s)}{\sststile{}{\alpha}\Gamma,\exists x\,A(x)}\quad \alpha_0<\alpha  \\  
& (\text{Cut}) \quad \prftree{\sststile{}{\alpha_0}\Gamma,A\quad \sststile{}{\alpha_0}\Gamma,\neg A}{\sststile{}{\alpha}\Gamma} \quad \alpha_0<\alpha\\
& (\Sigma\text{-Ref})\quad \prftree{\sststile{}{\alpha_0}\Gamma,A}{\sststile{}{\alpha}\Gamma,\exists z\,A^z} \quad \alpha_0,\Omega<\alpha,\,A\text{ is a $\Sigma$-formula}
\end{align*}
Proof trees that follow the rules of $\RS_\Omega(\bV)$ are called \textbf{$\RS_\Omega(\bV)$-proofs} or \textbf{$\RS_\Omega(\bV)$-derivations}. The \textbf{end sequent} of a $\RS_\Omega(\bV)$-derivation is the lowest sequent (i.e. conclusion) of that derivation. The \textbf{direct subderivation(s)} of a derivation are the derivation(s) from which the end sequent is inferred, if they exist. The \textbf{side formula(s)} of an inference are the formula(s) that are irrelevant to that inference. The \textbf{minor formula} of an inference is the formula to which the inference rule is applied, and the \textbf{principal formula} of that inference is the formula formed by that inference. \par 
For example, in the inference 
\begin{align*}
\prftree[r]{($\wedge$)}
{\Gamma,\phi\quad\Gamma,\psi}
{\Gamma,\phi\wedge \psi}
\end{align*}
the side formulas are $\Gamma$, the minor formulas are $\phi,\psi$ and the principal formula is $\phi\wedge \psi$. If a derivation $D$ ends with this inference, then the end sequent of $D$ is $\Gamma,\phi\wedge \psi$, and the direct subderivations of $D$ are the derivations which derived $\Gamma,\phi$ and $\Gamma,\psi$ respectively. \par 
The \textbf{rank} of a term or formula is defined as follows. 
\begin{itemize}
\item $\rank(u)=\omega\cdot |u|$. 
\item $\rank(\pm u\in v)=\max(\rank(u),\rank(v))+1$. 
\item $\rank(A\wedge B)=\rank(A\vee B)=\max(\rank(A),\rank(B))+1$. 
\item $\rank((\exists x\in u)F(x))=\rank((\forall x\in u)F(x) )=\max(\rank(u)+3,\rank(F(c_\emptyset))+2)$. 
\item $\rank(\exists x\,F(x))=\rank(\forall x \,F(x))=\max(\Omega,\rank(F(c_\emptyset))+1)$. 
\end{itemize}
$\sststile{\rho}{\alpha}\Gamma$ will be used to denote that $\sststile{}{\alpha}\Gamma$ and all cut formulas appearing in the derivation have rank $<\rho$. 
\end{named}

\begin{obsdep}
\label{rank observation}
For each formula $A$, define 
\begin{align*}
k(A)=\{ |t|\mid t\text{ occurs in }A\}\cup\{ \Omega\mid \text{if $A$ contains an unbounded quantifier}\}.
\end{align*}
\begin{enumerate}
\item For each formula $A$, $\rank(A)=\omega\cdot\max(k(A))+n$ for some $n<\omega$. 
\item $\rank(A)<\Omega$ iff $A$ is $\Delta_0$; thus $A$ has rank $\Omega$ iff $A$ is $\exists x\,F(x)$ or $\forall x\,F(x)$ where $F\in \Delta_0$. 
\end{enumerate}
\end{obsdep}

\begin{lem}
For each formula $A(s)$, if $|s|<\max(k(A(s)))$, then $\rank(A(s))=\rank(A(c_\emptyset))$. 
\end{lem}
\begin{proof}
By induction on complexity of $A$. 
\end{proof}

Some formulas can be regarded as generalisations of disjunctions or conjunctions:
\begin{itemize}
\item $A_0\wedge A_1\simeq \bigwedge_{i\in \{ 0,1\}}A_i$. 
\item $A_0\vee A_1\simeq \bigvee_{i\in \{0,1\}}A_i$. 
\item $(\forall x\in t)A(x)\simeq \bigwedge_s(s\in t\to A(s))$.
\item $(\exists x\in t)A(x)\simeq \bigvee_s(s\in t\wedge A(s))$. 
\item $\forall x\,A(x)\simeq \bigwedge_sA(s)$. 
\item $\exists x\,A(x)\simeq \bigvee_sA(s)$. 
\end{itemize}

\begin{lem}
\label{formula rank higher than subformulas}
If $A$ contains an unbounded quantifier and $A\simeq \bigvee_{i\in y} A_i$ or $A\simeq \bigwedge_{i\in y} A_i$ then 
\begin{align*}
(\forall i\in y)\rank(A_i)<\rank(A).
\end{align*}
\end{lem}
\begin{proof}
Straightforward by using the previous lemma. 
\end{proof}

The following lemmas will be formalised in Section \ref{sec:Formalisation of cut elimination in KP+TI}; their proofs outside KP can be carried out similarly to the standard ones (e.g. as in \cite{Schwichtenberg1977}). The parts that need consideration are how to code the infinitary proofs in KP and how to express and prove these lemmas in KP.

\begin{named}[Weakening]
\label{unformalised weakening}
If $\alpha\leq \alpha'$, $\rho\leq \rho'$, $\Gamma$ and $\Gamma'$ are finite sets of formulas, and $\sststile{\rho}{\alpha}\Gamma$, then $\sststile{\rho'}{\alpha'}\Gamma,\Gamma'$.
\end{named}

\begin{named}[Inversion]
If $A$ is not $\Delta_0$, $A\simeq\bigwedge (A_i)_{i\in y}$ and $\sststile{\rho}{\alpha}\Gamma,A$ then for all $i\in y$, $\sststile{\rho}{\alpha}\Gamma,A_i$. 
\end{named}

\begin{named}[Reduction]
Suppose $\rank(C)=\rho>\Omega$. If $\sststile{\rho}{\alpha}\Gamma,\neg C$ and $\sststile{\rho}{\beta}\Gamma,C$, then $\sststile{\rho}{\alpha+\beta}\Gamma$. 
\end{named}

\begin{named}[Cut elimination]
\label{unformalised cut elimination}
If $\sststile{\rho+1}{\beta}\Gamma$ and $\rho>\Omega$, then $\sststile{\rho}{\omega^\beta}\Gamma$. 
\end{named}

\section{Partial recursive set functions}

Most naively (and so not quite possibly), one might attempt to define the $\RSderivation$ in such a way: if for all terms $s$, $D_s\in \RSderivation$ with end sequent $\Gamma,F(s)$, then $\langle \code{\forall}, \{ D_s\mid s\text{ a term}\}\rangle$ is an $\RSderivation$ with end sequent $\Gamma,\forall x\,F(x)$. At least two problems are evident here: (1) $\{ D_s\mid s\text{ a term}\}$ is a proper class, so we need to enclose this class-many information into a set function expressed as an index; (2) this recursive definition involves an unbounded universal quantifier, which is not allowed in KP in general. \par

To fix these problems altogether, we use partial recursive set functions, for which the recursion theorem applies. Here we cite \cite{Rathjen2012} Section 2.2 to define partial $E$-recursive functions; more on partial $E$-recursive functions can be found in \cite{Normann1978} and \cite{Sacks2017} Chapter X.

\begin{defn}
Let $\mathbf{k,s,p,p_0,p_1,s_N,p_N,d_N,\bar{0},\bar{\pmb{\omega}},\pmb{\gamma},\pmb{\rho},\pmb{\nu},\pmb{\pi}},\mathbf{i_1,i_2,i_3}$ be distinct natural numbers; they will be indices for initial partial $E$-recursive functions. \par 
We define a class $\mathbb{E}$ of triples $\langle e,x,y\rangle$ by induction. Instead of $\langle e,x,y\rangle\in \mathbb{E}$ we will write $[e](x)\simeq y$; if $n>1$, we use $[e](x_1,\ldots,x_n)\simeq y$ to convey that 
\begin{align*}
[e](x_1)\simeq \langle e,x_1\rangle\wedge [\langle e,x_1\rangle](x_2)\simeq \langle e,x_1,x_2\rangle\wedge \ldots\wedge [\langle e,x_1,\ldots,x_{n-1}\rangle](x_n)\simeq y.
\end{align*}
We say that $[e](x)$ is defined, written $\convg{[e](x)}$, if $[e](x)\simeq y$ for some $y$. Let $\mathbb{N}$ denote $\omega$; $\mathbb{E}$ is defined by the following clauses:
\begin{itemize}
\item $[\mathbf{k}](x,y)\simeq x$. 
\item $[\mathbf{s}](x,y,z)\simeq [[x](z)]([y](z))$.\par 
$[\mathbf{s}](x,y,z)$ is not defined unless $[x](z),[y](z)$ and $[[x](z)]([y](z))$ are already defined; the clause for $\mathbf{s}$ should be read as a conjunction of the following clauses: $[\mathbf{s}](x)\simeq \langle \mathbf{s},x\rangle$, $[\langle \mathbf{s},x\rangle](y)\simeq \langle \mathbf{s},x,y\rangle$, and, if there exist $a,b,c$ such that $[x](z)\simeq a$, $[y](z)\simeq b$, $[a](b)\simeq c$, then $[\langle \mathbf{s},x,y\rangle](z)\simeq c$. 
\item $[\mathbf{p}](x,y)\simeq \langle x,y\rangle$. 
\item $[\mathbf{p_0}](x)\simeq (x)_0$. 
\item $[\mathbf{p_1}](x)\simeq (x)_1$. 
\item $[\mathbf{s_N}](n)\simeq n+1$ if $n\in \mathbb{N}$. 
\item $[\mathbf{p_N}](0)\simeq 0$. 
\item $[\mathbf{p_N}](n+1)\simeq n$ if $n\in \mathbb{N}$. 
\item $[\mathbf{d_N}](n,m,x,y)\simeq x$ if $n,m\in \mathbb{N}$ and $n=m$. 
\item $[\mathbf{d_N}](n,m,x,y)\simeq y$ if $n,m\in \mathbb{N}$ and $n\neq m$. 
\item $[\mathbf{\bar{0}}](x)\simeq 0$. 
\item $[\pmb{\bar{\omega}}](x)\simeq \omega$. 
\item $[\pmb{\pi}](x,y)\simeq \{ x,y\}$. 
\item $[\pmb{\nu}](x)\simeq \bigcup x$. 
\item $[\pmb{\gamma}](x,y)\simeq x\cap (\bigcap y)$. 
\item $[\pmb{\rho}](x,y)\simeq \{ [x](u)\mid u\in y\}$ if $[x](u)$ is defined for all $u\in y$. \par 
Similarly to the clause for $\mathbf{s}$, this means that $[\pmb{\rho}](x)\simeq \langle \pmb{\rho},x\rangle$, and if there is a function $f$ with domain $y$ such that $[x](u)\simeq f(u)$ for all $u\in y$, then $[\langle \pmb{\rho},x\rangle](y)\simeq \{ f(u)\mid u\in y\}$. 
\item $[\mathbf{i_1}](x,y,z)\simeq \{ u\in x\mid y\in z\}$. 
\item $[\mathbf{i_2}](x,y,z)\simeq \{ u\in x\mid u\in y\to u\in z\}$. 
\item $[\mathbf{i_3}](x,y,z)\simeq \{ u\in x\mid u\in y\to z\in u\}$. 
\end{itemize}
\end{defn}

\begin{prop}
\label{E can be defined}
The class $\mathbb{E}$ can be defined in KP; $\mathbb{E}$ is a $\Sigma$ class, and for all $e,x,y,y'$,
\begin{align*}
\langle e,x,y\rangle\in \mathbb{E}\wedge \langle e,x,y'\rangle\in \mathbb{E}\Rightarrow y=y'.
\end{align*}
\end{prop}
\begin{proof}
See \cite{Rathjen2012}; to prove the definability of $\mathbb{E}$, a theorem (\cite{AczelRathjen2001} Theorem 11.14) involving \emph{inductive definition} is quoted: the definition of $\mathbb{E}$ is an inductive definition by $\Sigma$ clauses, so the definability of $\mathbb{E}$ in KP follows from that theorem.
\end{proof}

\begin{defn}
\textbf{Application terms} are defined inductively as follows.
\begin{itemize}
\item The constants $\mathbf{k,s,p,p_0,p_1,s_N,p_N,d_N,\bar{0},\bar{\pmb{\omega}},\pmb{\gamma},\pmb{\rho},\pmb{\nu},\pmb{\pi}},\mathbf{i_1,i_2,i_3}$ are application terms. 
\item Variables are application terms. 
\item If $s$ and $t$ are application terms, then $(st)$ is an application term. 
\end{itemize}
A \textbf{closed} application term is an application term that does not contain variables.
\end{defn}

\begin{defn}
If $r$ is an application term and $u$ is a variable, we define the formula $[r\simeq u]^\wedge$ inductively as follows.
\begin{itemize}
\item If $r$ is a constant or a variable, $[r\simeq u]^\wedge$ is $r=u$. 
\item If $r$ is $(st)$, then $[r\simeq u]^\wedge$ is $\exists xy([s\simeq x]^\wedge\wedge [t\simeq y]^\wedge\wedge[x](y)\simeq u)$. 
\end{itemize}
\end{defn}

\begin{namedrmk*}[Notation]
$ $\par 
\begin{itemize}
\item $\convg{t}$ denotes $\exists x[t\simeq x]^\wedge$ (i.e. $t$ is defined).
\item $t(a_1,\ldots,a_n)\simeq b$ denotes
\begin{align*}
\exists x_1\ldots x_n\exists y(x_1=a_1\wedge \ldots\wedge x_n=a_n\wedge y=b\wedge [t(x_1,\ldots,x_n)\simeq y]^\wedge ). 
\end{align*}
\item $st_1\ldots t_n$ denotes $((\ldots(st_1)\ldots)t_n )$ (i.e. it is the functional $s(t_1,\ldots,t_n)$). 
\end{itemize}
\end{namedrmk*}

%
%

\begin{defn}
A partial $n$-place (class) function $\Upsilon$ is said to be a \textbf{partial $E$-recursive function} if there exists a closed application term $t_\Upsilon$ such that 
\begin{align*}
\dom(\Upsilon)=\{ (a_1,\ldots,a_n)\mid \convg{t_\Upsilon(a_1,\ldots,a_n)}\}
\end{align*}
and for all sets $(a_1,\ldots,a_n)\in \dom(\Upsilon)$, 
\begin{align*}
t_\Upsilon(a_1,\ldots,a_n)\simeq \Upsilon(a_1,\ldots,a_n). 
\end{align*}
In this case, $t_\Upsilon$ is said to be an \textbf{index} for $\Upsilon$. 
\end{defn}

If $\Upsilon_1,\Upsilon_2$ are partial $E$-recursive functions, then write $\Upsilon_1(\vec{a})\simeq \Upsilon_2(\vec{a})$ if neither $\Upsilon_1(\vec{a})$ nor $\Upsilon_2(\vec{a})$ are defined, or $\Upsilon_1(\vec{a})$ and $\Upsilon_2(\vec{a})$ are both defined and equal.

\begin{namedrmk}[Remark]
\label{partial recursive universal computation function}
Observe that, by using $[\mathbf{s}]$ (as well as pairing and projection functions), we have a partial recursive universal computation function $U$ that outputs $[e](x)$ on input $e,x$.
\end{namedrmk}

In the following lemma, the application term $t[x]$ means the application term $t$ with term $x$ indicated; $t[y]$ is the result of replacing every occurrence of $x$ in $t$ by $y$.

\begin{named}[Abstraction lemma]
For every application term $t[x]$ there exists an application term $\lambda x.t[x]$ such that the following holds:
\begin{align*}
\forall x_1\ldots x_n(\convg{\lambda x.t[x]}\wedge \forall y(\lambda x.t[x])y\simeq t[y]  ).
\end{align*}
\end{named}
\begin{proof}
By induction on the term $t$. If $t$ is the variable $x$, $\lambda x.x$ is $\mathbf{skk}$. If $t$ is a constant or a variable other than $x$, $\lambda x.t$ is $\mathbf{k} t$. If $t\equiv (uv)$, $\lambda x.uv$ is $\mathbf{s}(\lambda x.u )(\lambda x.v)$.
\end{proof}

\begin{named}[Recursion theorem]
There exists a closed application term $\mathsf{R}$ such that for any application terms $f,x$, 
\begin{align*}
\convg{\mathsf{R}f}\wedge \mathsf{R}fx\simeq f(\mathsf{R}f)x. 
\end{align*}
\end{named}
\begin{proof}
Take $\mathsf{R}$ to be $\lambda f.tt$, where $t$ is $\lambda y\lambda x.f(yy)x$. 
\end{proof}

In particular, given any partial $E$-recursive $f$, there is an index $e$ such that 
\begin{align*}
\Phi_e(x)\simeq f(e,x)
\end{align*}
for all $x$ ($\Phi_e$ stands for the same thing as $[e]$). \par

\begin{namedrmk}[Remark]
\label{partial E-recursive include primitive-recursive}
The reader may check that the class of partial $E$-recursive functions include the class of primitive-recursive set functions, defined in \ref{primitive-recursive set function}.
\end{namedrmk}


\begin{namedrmk}[Remark]
\label{partial E-recursive include self-reference}
We will describe a partial $E$-recursive function $f$ taking parameters $e$ and $x$, where $e$ is an index of a yet unknown function. By recursion theorem, we obtain an $e$ such that $\Phi_e$ behaves in the same way as $f(e,\cdot)$; thus, if we want to define a partial $E$-recursive function $g$, we may assume that $g$ already knows its own index. \par 
For example, if we want to find a partial $E$-recursive function $g$ such that $g(x)\simeq\langle e,x\rangle$ where $e$ is an index for $g$, we just define $f(e,x)=\langle e,x\rangle$, and apply $\mathsf{R}$ to $f$ to obtain an index $e$ such that $\Phi_e(x)\simeq f(e,x)\simeq \langle e,x\rangle$. Therefore we may define $g$ by simply saying ``$g(x)=\langle e,x\rangle$ where $e$ is an index of $g$''. 
\end{namedrmk}

In the sequel, we will just say ``partial recursive'' instead of ``partial $E$-recursive''.

\section{Expressibility of the infinitary derivations}

In this section, we show how to code the $\RS_\Omega(\bV)$-derivations within KP. \fullref{def of RS(V)-derivation outside KP} gives an idea of how the codes are structured, and the formal definition inside KP is given in \ref{def of RS(V)-derivation inside KP}.

\begin{defn}
\label{def of RS(V)-derivation outside KP}
Let us fix a natural coding function $\code{\cdot}$ for formulas and finite sets of formulas, as well as a set of symbols $\code{\mathrm{Axiom}}$, $\code{\wedge}$, $\code{\vee}$ etc.; we give the definition of codes of $\RS_\Omega(\bV)$-derivations \emph{outside} KP; we will be able to read off primitive-recursively from such a code $u$
\begin{enumerate}
\item the name of the last inference of the proof, 
\item its principal formula and side formulas, 
\item the end sequent,
\item a bound for the length of the proof,
\item a bound for the cut rank of the proof.  
\end{enumerate}
The corresponding primitive-recursive functions will be denoted by Rule$(u)$, PF$(u)$, SF$(u)$, End$(u)$, $\mathrm{Length}(u)$, Rank$(u)$, respectively. The ordinal notation system we use is the one $\langle R,\prec\rangle$ we defined in \ref{sec:Preliminary definitions}.
\begin{itemize}
\item If $\Gamma$ contains a true $\Delta_0$-formula, then
\begin{align}
\label{code Axiom}
\langle \code{\text{Axiom}},\code{\Gamma}\rangle   
\end{align}
is a member of $\RSderivation$. 
\item If $u,v\in \RSderivation$, $\mathrm{End}(u)=\code{\Gamma,A}$, $\mathrm{End}(v)=\code{\Gamma,B}$, and $\mathrm{Length}(u),\mathrm{Length}(v)\prec a$, $\max(\mathrm{Rank}(u),$ $\mathrm{Rank}(v))\preceq r$, then 
\begin{align}
\langle \code{\wedge},\code{A\wedge B},\code{\Gamma},a,r,u,v\rangle 
\end{align}
is a member of $\RSderivation$. 
\item If $u\in \RSderivation$, $\mathrm{End}(u)=\code{\Gamma,A}$, and $\mathrm{Length}(u)\prec a$, $\mathrm{Rank}(u)\preceq r$, then 
\begin{align}
\langle\code{\vee_0},\code{A\vee B},\code{\Gamma},a,r,u\rangle   
\end{align}
is a member of $\RSderivation$. 
\item If $u\in \RSderivation$, $\mathrm{End}(u)=\code{\Gamma,B}$, and $\mathrm{Length}(u)\prec a$, $\mathrm{Rank}(u)\preceq r$, then 
\begin{align}
\langle\code{\vee_1},\code{A\vee B},\code{\Gamma},a,r,u\rangle  
\end{align}
is a member of $\RSderivation$. 
\item Let $e$ be an index for a partial recursive set function. If for all $s$, $\Phi_e(s)=:u_s\in \RSderivation$, $\mathrm{End}(u_s)=\code{\Gamma,s\in t\to F(s)}$, $\mathrm{Length}(u_s)\prec a$, $\mathrm{Rank}(u_s)\preceq r$, then 
\begin{align}
\langle \code{b\forall},\code{(\forall x\in t)F(x)},\code{\Gamma},a,r,e\rangle 
\end{align}
is a member of $\RSderivation$. 
\item If $u\in \RSderivation$, $\mathrm{End}(u)=\code{\Gamma,s\in t\wedge F(s)}$ and $\mathrm{Length}(u)\prec a$, $\mathrm{Rank}(u)\preceq r$, then 
\begin{align}
\langle \code{b\exists}, \code{(\exists x\in t)F(x)}, \code{\Gamma}, a,r,u\rangle  
\end{align}
is a member of $\RSderivation$. 
\item Let $e$ be an index for a partial recursive set function. If for all $t$, $\Phi_e(t)=:u_t\in \RSderivation$, $\mathrm{End}(u_t)=\code{\Gamma,F(t)}$, $\mathrm{Length}(u_t)\prec a$, $\mathrm{Rank}(u_t)\preceq r$, then 
\begin{align}
\langle \code{\forall},\code{\forall x\,F(x)},\code{\Gamma},a,r,e\rangle   
\end{align}
is a member of $\RSderivation$. 
\item If $u\in \RSderivation$, $\mathrm{End}(u)=\code{\Gamma,F(s)}$, and $\mathrm{Length}(u)\prec a$, $\mathrm{Rank}(u)\preceq r$, then 
\begin{align}
\langle \code{\exists}, \code{\exists x\,F(x)},\code{\Gamma},a,r,u\rangle   
\end{align}
is a member of $\RSderivation$. 
\item If $u,v\in \RSderivation$, $\mathrm{End}(u)=\code{\Gamma,A}$, $\mathrm{End}(v)=\code{\Gamma,\neg A}$, $\mathrm{Length}(u),\mathrm{Length}(v)\prec a$, $\max(\rank(A)+1,\mathrm{Rank}(u),\mathrm{Rank}(v))\preceq r$, then 
\begin{align}
\label{code Cut}
\langle \code{\mathrm{Cut}}, \code{A}, \code{\Gamma}, a ,r,u,v\rangle   
\end{align}
is a member of $\RSderivation$. In this case, contrary to others, $A$ is the cut formula, so it does not appear in the end sequent of \eqref{code Cut}. 
\item If $u\in \RSderivation$, $\mathrm{End}(u)=\code{\Gamma,A}$, $\max(\mathrm{Length}(u),\Omega)\prec a$, $\mathrm{Rank}(u)\preceq r$, then 
\begin{align}
\label{code Sigma-Ref}
\langle \code{\Sigma\text{-Ref}},\code{\exists z\,A^z}, \code{\Gamma}, a,r,u\rangle 
\end{align}
is a member of $\RSderivation$. 
\end{itemize}
It is easy to see how the information extraction functions (End, Rank, Length, etc.) can be defined, regardless of such derivation codes being definable or not. 
\end{defn}

For the moment we haven't seen if KP can express $u\in \RSderivation$, but KP surely can express if some set looks like a member of $\RSderivation$ (i.e. in one of the forms \eqref{code Axiom}--\eqref{code Sigma-Ref}). \par

The following definitions are made in KP.

\begin{defn}
We say that $w$ is a \textbf{quasicode of a $\RS_\Omega(\bV)$-derivation}, $w\in \RSquasicode$, if $(w)_0\in \{\code{\text{Axiom}},\code{\wedge},\code{\vee_0},\code{\vee_1},\code{b\forall},\code{b\exists},\code{\forall},\code{\exists},\code{\text{Cut}}, \code{\Sigma\text{-Ref}}\}$, and 
\begin{itemize}
\item if $(w)_0=\code{\text{Axiom}}$, then $w=\langle \code{\text{Axiom}},\code{\Gamma}\rangle$ and $\Gamma$ is a finite set of formulas;
\item if $(w)_0=\code{\wedge}$, then $w=\langle \code{\wedge},\code{A\wedge B},\code{\Gamma},a,r,u,v\rangle$, where $A,B$ are formulas, $\Gamma$ is a finite set of formulas, $a,r$ are in our ordinal notation, and $u,v$ are sets;
\item if $(w)_0=\code{\vee_0}$, then $w=\langle \code{\vee_0},\code{A\vee B},\code{\Gamma},a,r,u\rangle$, where $A,B$ are formulas, $\Gamma$ is a finite set of formulas, $a,r$ are in our ordinal notation, and $u$ is a set;
\item if $(w)_0=\code{\vee_1}$, then $w=\langle \code{\vee_1},\code{A\vee B},\code{\Gamma},a,r,u\rangle$, where $A,B$ are formulas, $\Gamma$ is a finite set of formulas, $a,r$ are in our ordinal notation, and $u$ is a set;
\item if $(w)_0=\code{b\forall}$, then $w=\langle \code{b\forall},\code{(\forall x\in t)F(x)},\code{\Gamma},a,r,e\rangle$, where $t$ is a set, $F$ is a formula, $\Gamma$ is a finite set of formulas, $a,r$ are in our ordinal notation, and $e$ is a set;
\item if $(w)_0=\code{b\exists}$, then $w=\langle \code{b\exists},\code{(\exists x\in t)F(x)},\code{\Gamma},a,r,u\rangle$, where $t$ is a set, $F$ is a formula, $\Gamma$ is a finite set of formulas, $a,r$ are in our ordinal notation, and $u$ is a set;
\item if $(w)_0=\code{\forall}$, then $w=\langle \code{\forall},\code{\forall x\,F(x)},\code{\Gamma},a,r,e\rangle$, where $F$ is a formula, $\Gamma$ is a finite set of formulas, $a,r$ are in our ordinal notation, and $e$ is a set;
\item if $(w)_0=\code{\exists}$, then $w=\langle \code{\exists},\code{\exists x\,F(x)},\code{\Gamma},a,r,u\rangle$, where $F$ is a formula, $\Gamma$ is a finite set of formulas, $a,r$ are in our ordinal notation, and $u$ is a set;
\item if $(w)_0=\code{\text{Cut}}$, then $w=\langle \code{\text{Cut}},\code{A},\code{\Gamma},a,r,u,v\rangle$, where $A$ is a formula, $\Gamma$ is a finite set of formulas, $a,r$ are in our ordinal notation, and $u$ is a set;
\item if $(w)_0=\code{\Sigma\text{-Ref}}$, then $w=\langle \code{\Sigma\text{-Ref}},\code{\exists z\,A^z},\code{\Gamma},a,r,u\rangle$, where $A$ is a formula, $\Gamma$ is a finite set of formulas, $a,r$ are in our ordinal notation, and $u$ is a set. 
\end{itemize}
\end{defn}

It is primitive-recursive to determine if $w\in \RSquasicode$.\par 

If $w$ is a quasicode, we can still read off the information Rule$(w)$, PF$(w)$, etc.; this definition exists only to eliminate unnecessarily fussy case distinctions, e.g., if we know that $w\in \RSquasicode$ and $w=\langle \code{\wedge},\code{A\wedge B},\code{\Gamma},a,r,u,v\rangle$, then $\Gamma$ must be a finite set of formulas.  \par

An $\RS_\Omega(\bV)$-derivation figure is a well-founded tree. This means that we can label its nodes by finite sequences of sets. 

\begin{namedrmk*}[Notation]
Let $\FinSeq$ denote the class of all finite sequences of sets. If $\sigma\in \FinSeq$ and $x$ is a set, let $\sigma x$ abbreviate $\sigma^\frown\langle x\rangle$. 
\end{namedrmk*}

It is easy to define a partial recursive set function $N$ such that, if $u$ codes an $\RS_\Omega(\bV)$-derivation, $\sigma\in \FinSeq$, and $N(u,\sigma)$ is a subderivation of $u$, then for appropriate collections of sets $x$ (depending on the rule by which $N(u,\sigma)$ is derived), $N(u,\sigma x)$ are the direct subderivations of $N(u,\sigma)$. \par 


\begin{defn}
We define the partial recursive set function $N(w,\sigma)$, for $w\in \RSquasicode$ and $\sigma\in \FinSeq$. 
\begin{itemize}
\item $N(w,\langle\rangle)=w$. 
\item Suppose $(N(w,\sigma))_0\in \{\code{\vee_0}, \code{\vee_1},  \code{b\exists}, \code{\exists}, \code{\Sigma\text{-Ref}} \}$ and $(N(w,\sigma))_5=u$. Then $N(w,\sigma0)=u$. 
\item Suppose $(N(w,\sigma))_0\in \{\code{\wedge},\code{\text{Cut}} \}$ and $(N(w,\sigma))_5=u$, $(N(w,\sigma))_6=v$. Then $N(w,\sigma0)=u$ and $N(w,\sigma1)=v$.
\item Suppose $(N(w,\sigma))_0\in \{ \code{b\forall}, \code{\forall}\}$ and $(N(w,\sigma))_5=e$. Then $N(w,\sigma x)=\Phi_e(x)$ for all sets $x$ (see \fullref{partial recursive universal computation function}). 
\end{itemize}
In any other cases we set $N(w,\sigma)=\emptyset$. 
\end{defn}

\begin{namedrmk*}[Notation]
From now on, we denote $N(w,\sigma)$ by $w_\sigma$. 
\end{namedrmk*}

Before we proceed to give the formal definition of the class of $\RS_\Omega(\bV)$-derivations, we mention that one can define truth predicates for formulas with bounded complexity.

\begin{named}[Definition of the truth predicates]
$\True_{\Gamma}$, $\Gamma$ being a syntactical complexity, is defined primitive-recursively. 
\begin{itemize}
\item $\True_{\Delta_0}(\code{a\in b})$ iff $a\in b$ holds; the rest of the cases are $\wedge,\vee,\neg,(\forall x\in t), (\exists x\in t)$ defined similarly to below.  
\item $\True_{\Sigma_n}(\code{\phi})$ iff one of the following holds. 
\begin{itemize}
\item For some $\diamond\in \{ \wedge,\vee\}$ and $\psi,\theta\in \Sigma_n$, $\code{\phi}=\code{\psi\diamond\theta}$, and $\True_{\Sigma_n}(\code{\psi})\diamond \True_{\Sigma_n}(\code{\theta})$. 
\item For some $\psi\in \Pi_n$, $\code{\phi}=\code{\neg\psi}$, and $\neg\True_{\Pi_n}(\code{\psi})$. 
\item For some $\psi\in \Sigma_n$, $\code{\phi}=\code{(\exists x\in t)\psi(x)}$, and $(\exists x\in t)\True_{\Sigma_n}(\code{\psi(\dot{x})})$.
\item For some $\psi\in \Sigma_n$, $\code{\exists x\,\psi(x)}$, and $\exists x\,\True_{\Sigma_n}(\code{\psi(\dot{x})})$. 
\item $\True_{\Pi_{n-1}}(\code{\phi})$. 
\end{itemize}
\item $\True_{\Pi_n}(\code{\phi})$ is defined in a symmetrical way. 
\end{itemize}
\end{named}

Then in KP we have that for any formula $\phi$ with complexity $\Gamma$,
\begin{align*}
\mathrm{True}_\Gamma(\code{\phi})\leftrightarrow \phi.
\end{align*}

\begin{named}[Definition of an $\RS_\Omega(\bV)$-derivation within $\KP$]
\label{def of RS(V)-derivation inside KP}
A set $w$ is a code of a $\RS_\Omega(\bV)$-derivation, $w\in \RSderivation$, iff for all $\sigma\in \FinSeq$, the following holds:
\begin{itemize}
\item If $w_\sigma\neq \emptyset$ then $w_\sigma\in \RSquasicode$. 
\item If $w_\sigma=\langle \code{\text{Axiom}},\code{\Gamma}\rangle$, then $\Gamma$ contains a true $\Delta_0$-formula (i.e. there is $A\in \Gamma$ with complexity $\Delta_0$, and $\True_{\Delta_0}(\code{A})$); $w_{\sigma\tau}=\emptyset$ for any $\tau\neq \emptyset$.
\item If $w_\sigma=\langle \code{\wedge},\code{A\wedge B},\code{\Gamma},a,r,u,v\rangle$, then $\mathrm{End}(w_{\sigma0})=\code{\Gamma,A}$, $\mathrm{End}(w_{\sigma1})=\code{\Gamma,B}$, $\mathrm{Length}(w_{\sigma i})\prec a$ and $\mathrm{Rank}(w_{\sigma i})\preceq r$ for $i=0,1$; $w_{\sigma x}=\emptyset$ for $x\neq 0,1$. 
\item If $w_\sigma=\langle\code{\vee_0},\code{A\vee B},\code{\Gamma},a,r,u\rangle$, then $\mathrm{End}(w_{\sigma0})=\code{\Gamma,A}$, $\mathrm{Length}(w_{\sigma 0})\prec a$ and $\mathrm{Rank}(w_{\sigma 0})\preceq r$; $w_{\sigma x}=\emptyset$ for $x\neq 0$. 
\item If $w_\sigma=\langle\code{\vee_1},\code{A\vee B},\code{\Gamma},a,r,u\rangle$, then $\mathrm{End}(w_{\sigma0})=\code{\Gamma,B}$, $\mathrm{Length}(w_{\sigma 0})\prec a$ and $\mathrm{Rank}(w_{\sigma 0})\preceq r$; $w_{\sigma x}=\emptyset$ for $x\neq 0$. 
\item If $w_\sigma=\langle \code{b\forall},\code{(\forall x\in t)F(x)},\code{\Gamma},a,r,e\rangle$, then for all $s$, $\mathrm{End}(w_{\sigma s})=\code{\Gamma,s\in t\to A(s)}$, $\mathrm{Length}(w_{\sigma s})\prec a$, $\mathrm{Rank}(w_{\sigma s})\preceq r$.
\item If $w_\sigma=\langle \code{b\exists}, \code{(\exists x\in t)F(x)}, \code{\Gamma}, a,r,u\rangle$, then $\mathrm{End}(w_{\sigma0})=\code{\Gamma,s\in t\wedge A(s)}$ for some $s$, $\mathrm{Length}(w_{\sigma 0})\prec a$ and $\mathrm{Rank}(w_{\sigma 0})\preceq r$; $w_{\sigma x}=\emptyset$ for $x\neq 0$. 
\item If $w_\sigma=\langle \code{\forall},\code{\forall x\,F(x)},\code{\Gamma},a,r,e\rangle$, then for all $x$, $\mathrm{End}(w_{\sigma x})=\code{\Gamma,F(x)}$, $\mathrm{Length}(w_{\sigma x})\prec a$, $\mathrm{Rank}(w_{\sigma x})\preceq r$.
\item If $w_\sigma=\langle \code{\exists}, \code{\exists x\,F(x)},\code{\Gamma},a,r,u\rangle$, then $\mathrm{End}(w_{\sigma0})=\code{\Gamma,F(s)}$ for some $s$, $\mathrm{Length}(w_{\sigma 0})\prec a$ and $\mathrm{Rank}(w_{\sigma 0})\preceq r$; $w_{\sigma x}=\emptyset$ for $x\neq 0$. 
\item If $w_\sigma=\langle \code{\mathrm{Cut}}, \code{A}, \code{\Gamma}, a ,r,u,v\rangle$, then $\mathrm{End}(w_{\sigma0})=\code{\Gamma,A}$, $\mathrm{End}(w_{\sigma1})=\code{\Gamma,\neg A}$, $\mathrm{Length}(w_{\sigma i})\prec a$ and $\mathrm{Rank}(w_{\sigma i}),$ $\rank(A)+1\preceq r$ for $i=0,1$; $w_{\sigma x}=\emptyset$ for $x\neq 0,1$. 
\item If $w_\sigma=\langle \code{\Sigma\text{-Ref}},\code{\exists z\,A^z}, \code{\Gamma}, a,r,u\rangle$, then $\mathrm{End}(w_{\sigma0})=\code{\Gamma,A}$, $\mathrm{Length}(w_{\sigma 0})\prec a$, $\mathrm{Rank}(w_{\sigma 0})\preceq r$; $w_{\sigma x}=\emptyset$ for $x\neq 0$. 
\end{itemize}
\end{named}

The reader may check that this is a $\Pi_2$ definition (recall that $\mathbb{E}$ is $\Sigma$ \ref{E can be defined}). 

\section{Formalisation of cut elimination in \texorpdfstring{$\KP+\TI(\epsilon_{\Omega+1})$}{KP+TI(epsilon\_Omega+1)}}

\label{sec:Formalisation of cut elimination in KP+TI}

The lemmas \ref{unformalised weakening}--\ref{unformalised cut elimination} are now going to be formalised in $\KP+\TI(\epsilon_{\Omega+1})$; we will define partial recursive set functions Wkn (for ``weakening''), Inv (for ``inversion''), Red (for ``reduction''), CutElim (for ``cut elimination'') that transform members of $\RS_\Omega(\bV)$-derivation into the required forms. We do not care what these functions do to the sets which are not members of $\RSderivation$; they may simply diverge. In defining these functions, we may assume that they know their own indices (by \ref{partial E-recursive include self-reference}).

The reader is advised to familiarise themselves with the standard cut elimination proofs before proceeding; though the reasoning is exactly the same, our proofs may have been obscured in formalisation.\par 

If there is no danger of confusion, we drop the coding notation $\code{\cdot}$. 

\begin{named}[Weakening]
$\KP+\TI(\epsilon_{\Omega+1})$ proves the following. There is a partial recursive set function $f$ taking parameter $\Gamma'$ such that, whenever $\Gamma'$ is a finite set of $\RS_\Omega(\bV)$-formulas, $w\in\RSderivation$, $\mathrm{End}(w)=\Gamma$, $\mathrm{Length}(w)=a\prec \epsilon_{\Omega+1}$, we have $f(\Gamma', w)\in \RSderivation$, $\mathrm{End}(f(\Gamma', w))=\Gamma,\Gamma'$, $\mathrm{Length}(w)\succeq \mathrm{Length}(f(\Gamma',w))$ and $\mathrm{Rank}(w)\succeq \mathrm{Rank}(f(\Gamma',w))$.  
\end{named}
The function $f$ will be referred to as Wkn. \par 
\begin{proof}
In the following we suppress the parameter $\Gamma'$ and write simply $f(w)$. \par 
If $w\notin \RSquasicode$, we set $f(w)=\emptyset$. Now we assume that $w\in \RSquasicode$. \par 
If $w=\langle \code{\text{Axiom}},\code{\Gamma}\rangle$, then\\
$f(w)=\langle \code{\text{Axiom}},\code{\Gamma,\Gamma'}\rangle$.\par 
If $w=\langle \code{\wedge},\code{A\wedge B},\code{\Gamma},a,r,u,v\rangle$ and the lengths of $u,v$ are both $\prec a$, then\\
$f(w)=\langle \code{\wedge},\code{A\wedge B},\code{\Gamma,\Gamma'},a,r,f(u),f(v)\rangle$ (using the recursion theorem).\par 
If $w=\langle\code{\vee_0},\code{A\vee B},\code{\Gamma},a,r,u\rangle$ and the length of $u$ is $\prec a$, then \\
$f(w)=\langle \code{\vee_0},\code{A\vee B},\code{\Gamma,\Gamma'},a,r,f(u)\rangle$. The case of $\vee_1$ is similar. \par 
If $w=\langle \code{b\forall},\code{(\forall x\in t)F(x)},\code{\Gamma},a,r,e\rangle$, then \\
$f(w)=\langle \code{b\forall},\code{(\forall x\in t)F(x)},\code{\Gamma,\Gamma'},a,r,e'\rangle$ where $e'$ is the natural index of $f\circ \Phi_e$. \par 
If $w=\langle \code{b\exists}, \code{(\exists x\in t)F(x)}, \code{\Gamma}, a,r,u\rangle$ and the length of $u$ is $\prec a$, then \\
$f(w)=\langle \code{b\exists}, \code{(\exists x\in t)F(x)}, \code{\Gamma,\Gamma'}, a,r,f(u)\rangle$. \par 
If $w=\langle \code{\forall},\code{\forall x\,F(x)},\code{\Gamma},a,r,e\rangle$, then\\
$f(w)=\langle \code{\forall},\code{\forall x\,F(x)},\code{\Gamma,\Gamma'},a,r,e'\rangle$ where $e'$ is the natural index of $f\circ \Phi_e$. \par 
If $w=\langle \code{\exists}, \code{\exists x\,F(x)},\code{\Gamma},a,r,u\rangle$ and the length of $u$ is $\prec a$, then \\
$f(w)=\langle \code{\exists}, \code{\exists x\,F(x)},\code{\Gamma,\Gamma'},a,r,f(u)\rangle$. \par 
If $w=\langle \code{\mathrm{Cut}}, \code{A}, \code{\Gamma}, a ,r,u,v\rangle$ and the lengths of $u,v$ are $\prec a$, then \\
$f(w)=\langle \code{\mathrm{Cut}}, \code{A}, \code{\Gamma,\Gamma'}, a ,r,f(u),f(v)\rangle$ if the cut formula $A$ is not in $\Gamma'$, and \\
$f(w)=f(u)$ if $A\in \Gamma'$ (in our definition we assume $A\in \mathrm{End}(u)$).\par 
If $w=\langle \code{\Sigma\text{-Ref}},\code{\exists z\,A^z}, \code{\Gamma}, a,r,u\rangle$ and the length of $u$ is $\prec a$, then \\
$f(w)=\langle \code{\Sigma\text{-Ref}},\code{\exists z\,A^z}, \code{\Gamma,\Gamma'}, a,r,f(u)\rangle$. \par 
One can then prove using $\TI(\epsilon_{\Omega+1})$ that $f$ satisfies the requirements. 
\end{proof}

Note that in the cases $b\forall$ and $\forall$, we do not check if say $\Phi_e(s)$ has length $\prec a$ for every set $s$. It is unnecessary, and we cannot do so recursively either. \par

\begin{named}[Inversion]
$\KP+\TI(\epsilon_{\Omega+1})$ proves the following. There is a partial recursive set function $g$ taking parameter $A$ such that, whenever $A$ is not $\Delta_0$, $A\simeq \bigwedge_{i\in y}A_i$, $w\in\RSderivation$, $\mathrm{End}(w)=\Gamma,A$, $\mathrm{Length}(w)=a\prec \epsilon_{\Omega+1}$, we have $g(A,w,i)\in \RSderivation$, $\mathrm{End}(g(A,w,i))=\Gamma,A_i$, $\mathrm{Length}(w)\succeq \mathrm{Length}(g(A,w,i))$ and $\mathrm{Rank}(w)\succeq \mathrm{Rank}(g(A,w,i))$, for all $i\in y$.
\end{named}

The function $g$ will be referred to as Inv. 
\begin{proof}
In the following we suppress the parameter $A\simeq \bigwedge_{i\in y}A_i$. \par 
If $w\notin \RSquasicode$ or $A\notin \mathrm{End}(w)$, then we set $g(w,i)=\emptyset$ for all $i$. Now we assume that $w\in \RSquasicode$ and $A\in \mathrm{End}(w)$. \par 
If $w=\langle \code{\text{Axiom}},\code{\Gamma,A}\rangle$, then\\
$g(w,i)=\langle \code{\text{Axiom}},\code{\Gamma,A_i}\rangle$ for all $i\in y$. \par 
Now suppose $w$ is not an axiom. We first assume that $A$ is not the principal formula of the last inference of $w$. \par 
If $w=\langle \code{\wedge},\code{B},\code{\Gamma,A},a,r,u,v\rangle$, and the lengths of $u,v$ are $\prec a$, then \\
$g(w,i)=\langle \code{\wedge},\code{B},\code{\Gamma,A_i},a,r,g(u,i),g(v,i)\rangle$. \par 
If $w=\langle \code{\vee_0},\code{B},\code{\Gamma,A},a,r,u\rangle$, and the length of $u$ is $\prec a$, then \\
$g(w,i)=\langle \code{\vee_0},\code{B},\code{\Gamma,A_i},a,r,g(u,i)\rangle$. The case of $\vee_1$ is similar. \par 
If $w=\langle \code{b\forall},\code{B},\code{\Gamma,A_i},a,r,e\rangle$, then\\
$g(w,i)=\langle\code{b\forall},\code{B},\code{\Gamma,A_i},a,r,e'\rangle$, where $e'$ is the natural index of $g(\Phi_e(\cdot),i)$. \par 
The rest of the cases are similar to above. \par 
Now assume that $A$ is the principal formula of the last inference of $w$. \par 
If $A\equiv A_0\wedge A_1$, $w=\langle \code{\wedge},\code{A},\code{\Gamma},a,r,u,v\rangle$, and the lengths of $u,v$ are $\prec a$, then \\
$g(w,0)=g(u,0)$ if $A\in \Gamma$, $g(w,0)=u$ if $A\notin \Gamma$, and similarly\\
$g(w,1)=g(v,1)$ if $A\in \Gamma$, $g(w,1)=v$ if $A\notin \Gamma$. \par 
If $A\equiv (\forall x\in t)F(x)$, $w=\langle \code{b\forall},\code{A},\code{\Gamma},a,r,e\rangle$, then, for each $s$, \\
$g(w,s)=g(\Phi_e(s),s)$ if $A\in \Gamma$ and $\Phi_e(s)$ has length $\prec a$, $g(w,s)=\Phi_e(s)$ if $A\notin \Gamma$. \par 
If $A\equiv \forall x\,F(x)$, $w=\langle \code{\forall}, \code{\forall x\,F(x)},\code{\Gamma},a,r,e\rangle$, then, for each $s$, \\
$g(w,s)=g(\Phi_e(s),s)$ if $A\in \Gamma$ and $\Phi_e(s)$ has length $\prec a$, $g(w,s)=\Phi_e(s)$ if $A\notin \Gamma$. \par 
Again using $\TI(\epsilon_{\Omega+1})$ one checks that $g$ satisfies the requirements for $w\in \RSderivation$ with length $\prec \epsilon_{\Omega+1}$. 
\end{proof}

Recall the definition of natural sum $\alpha\#\beta$ for ordinals \ref{natural sum definition}.

\begin{named}[Reduction]
$\KP+\TI(\epsilon_{\Omega+1})$ proves the following. There is a partial recursive set function $h$ taking parameter $C$ such that, whenever $\rank(C)=r\succ \Omega$, $w_0,w_1\in \RSderivation$, $\mathrm{End}(w_0)=\Gamma,C$, $\mathrm{End}(w_1)=\Gamma,\neg C$, $\mathrm{Length}(w_0)=a\prec \epsilon_{\Omega+1}$, $\mathrm{Length}(w_1)=b\prec \epsilon_{\Omega+1}$, $\mathrm{Rank}(w_0)$ and $\mathrm{Rank}(w_1)$ are both $\preceq r$, we have $h(C,w_0,w_1)\in \RSderivation$, $\mathrm{End}(h(C,w_0,w_1))=\Gamma$, $\mathrm{Length}(h(C,w_0,w_1))\preceq a\#b$, $\mathrm{Rank}(h(C,w_0,w_1))\preceq r$. 
\end{named}
The function $h$ will be referred to as Red.
\begin{proof}
We suppress the parameter $C$ with rank $r\succ \Omega$. We assume that $w_0,w_1\in \RSquasicode$ satisfy $\mathrm{End}(w_0)=\Gamma,C$, $\mathrm{End}(w_1)=\Gamma,\neg C$, $\mathrm{Length}(w_0)=a\prec \epsilon_{\Omega+1}$, $\mathrm{Length}(w_1)=b\prec \epsilon_{\Omega+1}$, $\mathrm{Rank}(w_0)$ and $\mathrm{Rank}(w_1)$ are both $\preceq r$. The function $h$ is defined by recursion on $a\#b$ along $\epsilon_{\Omega+1}$. \par 
If $(w_0)_0=\code{\text{Axiom}}$ or $(w_1)_0=\code{\text{Axiom}}$, then\\
$h(w_0,w_1)=\langle \code{\text{Axiom}}, \code{\Gamma}\rangle$. \par 
Now we assume that neither of $w_0,w_1$ is an axiom. \par 
First assume that $C$ is not the principal formula of the last inference of $w_0$. We go through the ideas of two representative cases. \par 
Suppose $w_0=\langle \code{\wedge},\code{A_0\wedge A_1},\code{\Gamma,C},a,r,u,v\rangle$; then, assuming this inference is correct, $\mathrm{End}(u)=A_0,\Gamma,C$ and $\mathrm{End}(v)=A_1,\Gamma,C$. By weakening, $\mathrm{Wkn}(A_i,w_1)$ has end sequent $A_i,\Gamma,\neg C$ for $i=0,1$, so if both $u,v$ have lengths $\prec a$, then we may apply $h$ to get $h(u,\mathrm{Wkn}(A_0,w_1))$ and $h(v,\mathrm{Wkn}(A_1,w_1))$ whose end sequents are $A_0,\Gamma$ and $A_1,\Gamma$ respectively. Thus we may apply $\wedge$ to $h(u,\mathrm{Wkn}(A_0,w_1))$ and $h(v,\mathrm{Wkn}(A_1,w_1))$ to obtain $h(w_0,w_1)$. \par 
Suppose $w_0=\langle \code{\forall},\code{\forall x\,F(x)},\code{\Gamma,C},a,r,e\rangle$. If this inference is correct, then $\Phi_e(s)$ has end sequent $F(s),\Gamma,C$ for every set $s$; also $\mathrm{Wkn}(F(s),w_1)$ has end sequent $F(s),\Gamma,\neg C$ for every $s$. Therefore if $\Phi_e(s)$ has length $\prec a$ then $h(\Phi_e(s),\mathrm{Wkn}(F(s),w_1))$ has end sequent $F(s),\Gamma$, thus we may apply $\forall$ to the derivations $h(\Phi_e(s),\mathrm{Wkn}(F(s),w_1))$ to obtain $h(w_0,w_1)$. \par 
Now we proceed to the formal definition. \par 
If $w_0=\langle \code{\wedge},\code{A_0\wedge A_1},\code{\Gamma,C},a,r,u,v\rangle$ and $u,v$ have lengths $\prec a$, then \\
$h(w_0,w_1)=\langle \code{\wedge},\code{A_0\wedge A_1},\code{\Gamma}, a\# b, r, h(u,\mathrm{Wkn}(A_0,w_1)), h(v,\mathrm{Wkn}(A_1,w_1))\rangle$. \par  
If $w_0=\langle \code{\vee_0},\code{A\vee B},\code{\Gamma,C},a,r,u\rangle$ and $u$ has length $\prec a$, then \\
$h(w_0,w_1)=\langle \code{\vee_0},\code{A\vee B},\code{\Gamma},a\#b,r,h(u,\mathrm{Wkn}(A,w_1))\rangle$. The case of $\vee_1$ is similar. \par 
If $w_0=\langle \code{b\forall},\code{(\forall x\in t)F(x)},\code{\Gamma,C},a,r,e\rangle$, then,\\
$h(w_0,w_1)=\langle \code{b\forall},\code{(\forall x\in t)F(x)},\code{\Gamma},a\#b,r,e'\rangle$, where $e'$ is the natural index of\\
$\lambda s. h(\Phi_e(s),\mathrm{Wkn}(s\in t\to F(s),w_1))$. \par 
The other cases, as well as all the cases in which $\neg C$ is not the principal formula of the last inference of $w_1$, can be similarly dealt with. \par 
For the remainder we assume that $C$ is the principal formula of the last inference of $w_0$ and $\neg C$ is the principal formula of the last inference of $w_1$. We go through the idea of one representative case. \par 
Suppose $w_0=\langle \code{\vee_0},\code{A_0\vee A_1},\code{\Gamma},a,r,u\rangle$. Then $\mathrm{Wkn}(C,u)$ has end sequent $A_0,\Gamma,C$; $\mathrm{Wkn}(A_0,w_1)$ has end sequent $A_0,\Gamma,\neg C$. Therefore if $u$ has length $\prec a$, we may apply $h$ and the derivation $h(\mathrm{Wkn}(C,u),\mathrm{Wkn}(A_0,w_1))$ has end sequent $A_0,\Gamma$ with length $\prec a\#b$. Now $\mathrm{Inv}(\neg C,w_1,0)$ has end sequent $\neg A_0,\Gamma$, so we may apply a cut to $h(\mathrm{Wkn}(C,u),\mathrm{Wkn}(A_0,w_1))$ and $\mathrm{Inv}(\neg C,w_1,0)$. \par 
The following is the formal definition. \par 
If $w_0=\langle \code{\vee_0},\code{A_0\vee A_1},\code{\Gamma},a,r,u\rangle$ and $u$ has length $\prec a$, then\\
$h(w_0,w_1)=\langle \code{\text{Cut}},\code{A_0},\code{\Gamma},a\#b,r,h(\mathrm{Wkn}(C,u),\mathrm{Wkn}(A_0,w_1)),\mathrm{Inv}(\neg C,w_1,0)\rangle$. \par 
If $w_0=\langle \code{b\exists},\code{(\exists x\in t)F(x)},\code{\Gamma},a,r,u\rangle$, $u$ has length $\prec a$ and end sequent $s\in t\wedge F(s),\Gamma$ for some set $s$, then \\
$h(w_0,w_1)=\langle \code{\text{Cut}},\code{s\in t\wedge F(s)},\code{\Gamma},a\#b,r,h(\mathrm{Wkn}(C,u),\mathrm{Wkn}(s\in t\wedge F(s),w_1)), \mathrm{Inv}(\neg C,w_1,s)\rangle$.\\
(Note that since $\rank(C)\succ \Omega$, $F$ contains an unbounded quantifier, so \ref{formula rank higher than subformulas} applies.) \par 
The other cases are similar or symmetrical to the ones above. 
\end{proof}

\begin{named}[Cut elimination]
$\KP+\TI(\epsilon_{\Omega+1})$ proves the following. There is a partial recursive set function $j$ such that, whenever $w\in \RSderivation$, $\mathrm{End}(w)=\Gamma$, $\mathrm{Length}(w)\preceq a\prec \epsilon_{\Omega+1}$, $\mathrm{Rank}(w)\preceq r+1$ with $r\succ \Omega$, we have $j(w)\in \RSderivation$, $\mathrm{End}(j(w))=\Gamma$, $\mathrm{Length}(j(w))\preceq \omega^a$, $\mathrm{Rank}(j(w))\preceq r$. 
\end{named}

The function $j$ here will be referred to as CutElim.
\begin{proof}
We assume $w\in \RSquasicode$, $\mathrm{End}(w)=\Gamma$, $\mathrm{Length}(w)\preceq a\prec \epsilon_{\Omega+1}$, $\mathrm{Rank}(w)\preceq r+1$ with $r\succ \Omega$. \par 
We give the definition of $j$ only for the representative cases.\par 
If $w=\langle \code{\text{Axiom}},\code{\Gamma}\rangle$, then $j(w)=w$. \par 
If $w=\langle \code{\wedge},\code{A_0\wedge A_1},\code{\Gamma},a,r+1,u,v\rangle$, and $u,v$ have lengths $\prec a$, then \\
$j(w)=\langle \code{\wedge},\code{A_0\wedge A_1},\code{\Gamma},\omega^a,r,j(u),j(v)\rangle$. The cases of $\vee_0,\vee_1,b\exists,\exists,\Sigma\text{-Ref}$ are similar.\par 
If $w=\langle \code{b\forall},\code{(\forall x\in t)F(x)},\code{\Gamma},a,r+1,e\rangle$, then \\
$j(w)=\langle \code{b\forall},\code{(\forall x\in t)F(x)},\code{\Gamma},\omega^a,r,e'\rangle$, where $e'$ is the natural index of $j\circ \Phi_e$. The case of $\forall$ is similar. \par 
If $w=\langle \code{\text{Cut}},\code{A},\code{\Gamma},a,r+1,u,v\rangle$ and $u,v$ have lengths $\prec a$, then \\
$j(w)=\langle \code{\text{Cut}},\code{A},\code{\Gamma},a,r,j(u),j(v)\rangle$ if $\rank(A)\prec r$, $j(w)=\mathrm{Red}(j(u),j(v))$ if $\rank(A)=r$. 
\end{proof}

\section{Embedding \texorpdfstring{$\KP$ into $\RS_\Omega(\bV)$}{KP into RS(V)}}
\label{sec:Embedding KP}


We are going to formalise, within KP, that KP can be embedded into $\RS_\Omega(\bV)$; roughly, there is a partial recursive set function $P$ such that, if $p$ is a KP-proof, then $P(p)$ is an $\RSderivation$ with the same conclusion (see \ref{embedding theorem}). Combining this with our formalisation of cut elimination, every KP-proof can be transformed into an $\RS_\Omega(\bV)$-proof with their cuts partially eliminated. Along these proof trees, we are able to carry out a transfinite induction along $\epsilon_{\Omega+1}$ which gives us the result $\KP+\TI(\epsilon_{\Omega+1})\vdash \RFN(\KP)$.

\begin{namedrmk*}[Notation]
$ $\par 
\begin{itemize}
\item If $A$ is any $\RS_\Omega(\bV)$-formula, then $\mathrm{no}(A)=\omega^{\rank(A)}$.
\item If $\Gamma=\{ A_1,\ldots,A_n\}$ is a set of $\RS_\Omega(\bV)$-formulas, then $\mathrm{no}(\Gamma):=\mathrm{no}(A_1)\#\ldots \#\mathrm{no}(A_n)$.
\item $\dststile{}{}\Gamma$ abbreviates $\sststile{0}{\mathrm{no}(\Gamma)}\Gamma$, and $\dststile{\rho}{\alpha}\Gamma$ abbreviates $\sststile{\rho}{\mathrm{no}(\Gamma)\#\alpha}\Gamma$. 
\end{itemize}
\end{namedrmk*}

%
%


\begin{named}[Embedding lemmas]
\label{embedding lemmas}
$ $\par 
\begin{enumerate}[label=\arabic*.]
\item $\dststile{}{}A,\neg A$.
\item (Extensionality). $\dststile{}{}s_1\neq t_1,\ldots,s_n\neq t_n,\neg A(s_1,\ldots,s_n),A(t_1,\ldots,t_n)$. 
\item (Set induction). $\dststile{}{\omega^{\rank(A)}}A\to \forall x\,F(x)$, where $A\equiv \forall x((\forall y\in x)F(y)\to F(x) )$. 
\item (Pair). $\dststile{}{}\exists z(s\in z\wedge t\in z)$. 
\item (Union). $\dststile{}{}\exists z(\forall y\in s)(\forall x\in y)(x\in z)$. 
\item (Infinity). $\dststile{}{}\exists x((\exists z\in x)z\in x\wedge (\forall y\in x)(\exists z\in x)y\in z )$. 
\item ($\Delta_0$-separation). If $A$ is $\Delta_0$, $\dststile{}{}\exists y((\forall x\in y)(x\in s\wedge A(x,\vec{t}))\wedge (\forall x\in s)(A(x,\vec{t})\to x\in y ) )$. 
\item ($\Delta_0$-collection). If $F$ is $\Delta_0$, $\dststile{}{}(\forall x\in s)\exists y\,F(x,y,\vec{t})\to \exists z(\forall x\in s)(\exists y\in z)F(x,y,\vec{t})$. 
\end{enumerate}
\end{named}

\begin{sublemma}[Lemma]
\label{LEM lemma}
KP proves the following. There is a partial recursive set function $f$ such that, for any $\RS_\Omega(\bV)$-formula $A$, $f(A)\in \RSderivation$, $\mathrm{End}(f(A))=A,\neg A$, $\mathrm{Length}(f(A))\preceq \mathrm{no}(A,\neg A)$, $\mathrm{Rank}(f(A))=0$. 
\end{sublemma}
The $f$ here will be referred to as LEM. 
\begin{proof}
We define $f$ by recursion on complexity of $A$. \par 
If $A$ is $\Delta_0$, we set \\
$f(A)=\langle \code{\text{Axiom}},\code{A,\neg A}\rangle$. From now on we assume that $A$ is not $\Delta_0$. \par 
If $A$ is $A_0\vee A_1$, let \\
$w_0=\langle \code{\vee_0},\code{A_0\vee A_1},\code{\neg A_0},\mathrm{no}(A_0,\neg A_0)+1,0,f(A_0)\rangle$,\\
$w_1=\langle \code{\vee_1},\code{A_0\vee A_1},\code{\neg A_1},\mathrm{no}(A_1,\neg A_1)+1,0,f(A_1)\rangle$, then set\\
$f(A)=\langle \code{\wedge},\code{\neg A_0\wedge \neg A_1},\code{A_0\vee A_1},\mathrm{no}(A,\neg A),0,w_0,w_1\rangle$. \par 
If $A$ is $(\exists x\in t)F(x)$, write $B(s)$ for $s\in t\wedge F(s)$, and let $e$ be the natural index of the function\\
$\Phi_e(s)=\langle \code{b\exists},\code{(\exists x\in t)F(x)},\code{s\in t\to \neg F(s)},\mathrm{no}(B(s),\neg B(s) )+1,0,f(B(s))\rangle$; then we set\\
$f(A)=\langle \code{b\forall},\code{(\forall x\in t)\neg F(x)},\code{(\exists x\in t)F(x)},\mathrm{no}(A,\neg A),0,e\rangle$. \par 
If $A$ is $\exists x\,F(x)$, let $e$ be the natural index of the function\\
$\Phi_e(s)=\langle \code{\exists},\code{\exists x\,F(x)},\code{\neg F(s)},\mathrm{no}(F(s),\neg F(s) )+1,0,f(F(s))\rangle$, and set\\
$f(A)=\langle \code{\forall},\code{\forall x\,\neg F(x)},\code{\exists x\,F(x)},\mathrm{no}(A,\neg A),0,e\rangle$. \par 
The other cases are symmetrical. 
\end{proof}


\begin{sublemma}[Extensionality]
KP proves the following. There is a partial recursive set function Ext such that, whenever $A(a_1,\ldots,a_n)$ is a KP-formula with free variables among $a_1,\ldots,a_n$, $\vec{s}:=\langle s_1,\ldots,s_n\rangle$, $\vec{t}:=\langle t_1,\ldots,t_n\rangle$ are sequences of terms with length $n$, we have $\mathrm{Ext}(A,\vec{s},\vec{t})\in \RSderivation$, $\mathrm{End}(\mathrm{Ext}(A,\vec{s},\vec{t}))=s_1\neq t_1,\ldots,s_n\neq t_n,\neg A(\vec{s}),A(\vec{t})$, $\mathrm{Length}(\mathrm{Ext}(A,\vec{s},\vec{t}))\preceq \mathrm{no}(s_1\neq t_1,\ldots,s_n\neq t_n,\neg A(\vec{s}),A(\vec{t}))$, $\mathrm{Rank}(\mathrm{Ext}(A,\vec{s},\vec{t}))=0$. 
\end{sublemma}
\begin{proof}
Define Ext by recursion on complexity of $A$; we will write $\vec{s}\neq \vec{t}$ for $s_1\neq t_1,\ldots,s_n\neq t_n$. \par 
If $A$ is $\Delta_0$, set $\mathrm{Ext}(A,\vec{s},\vec{t})=\langle \code{\text{Axiom}},\code{\vec{s}\neq \vec{t},\neg A(\vec{s}),A(\vec{t})}\rangle$. In this case, if $\mathrm{Ext}(A,\vec{s},\vec{t})\notin \RSderivation$, then by definition of $\RSderivation$, we have 
\begin{align*}
\neg\mathrm{True}_{\Delta_0}(\code{\vec{s}\neq \vec{t}})\wedge \neg\mathrm{True}_{\Delta_0}(\code{\neg A(\vec{s})})\wedge \neg\mathrm{True}_{\Delta_0}(\code{A(\vec{t})}),
\end{align*}
which implies $\vec{s}=\vec{t}\wedge A(\vec{s})\wedge \neg A(\vec{t})$ and yields a contradiction in KP. \par 
From now on we assume that $A$ is not $\Delta_0$. \par 
If $A$ is $A_0\wedge A_1$, $\mathrm{Ext}(A,\vec{s},\vec{t})$ is defined in the following way:
\begin{align*}
\prftree[l]{($\wedge$)}
{\prftree[l]{($\vee$)}
{\mathrm{Ext}(A_0,\vec{s},\vec{t})\dststile{}{}\vec{s}\neq \vec{t},\neg A_0(\vec{s}),A_0(\vec{t})}
{\vec{s}\neq \vec{t},\neg A_0(\vec{s})\vee \neg A_1(\vec{s}),A_0(\vec{t})}
\qquad
\prftree[l]{($\vee$)}
{\mathrm{Ext}(A_1,\vec{s},\vec{t})\dststile{}{}\vec{s}\neq \vec{t},\neg A_1(\vec{s}),A_1(\vec{t})}
{\vec{s}\neq \vec{t},\neg A_0(\vec{s})\vee \neg A_1(\vec{s}),A_1(\vec{t})}
}
{\vec{s}\neq\vec{t},\neg A_0(\vec{s})\vee \neg A_1(\vec{s}),A_0(\vec{t})\wedge A_1(\vec{t})}
\end{align*}
If $A$ is $(\forall x\in t)F(x)$, we derive $\vec{s}\neq \vec{t},(\exists x\in t(\vec{s}))\neg F(x,\vec{s}),s\notin t(\vec{t})\vee F(s,\vec{t})$ uniformly for all $s$:
\begin{align*}
\prftree[l]{($b\exists$)}
{\prftree[l]{($\wedge$)}
{\prftree[l]{($\vee$)}
{\mathrm{Ext}(s\in t(\vec{a}),\vec{s},\vec{t} )\dststile{}{}\vec{s}\neq \vec{t},s\in t(\vec{s}),s\notin t(\vec{t})}
{\vec{s}\neq \vec{t},s\in t(\vec{s}),s\notin t(\vec{t})\vee F(s,\vec{t})}
\qquad
\prftree[l]{($\vee$)}
{\mathrm{Ext}(F(s,\vec{a}),\vec{s},\vec{t})\dststile{}{}\vec{s}\neq \vec{t},\neg F(s,\vec{s}), F(s,\vec{t})}
{\vec{s}\neq \vec{t},\neg F(s,\vec{s}),s\notin t(\vec{t})\vee F(s,\vec{t})}
}
{\vec{s}\neq \vec{t}, s\in t(\vec{s})\wedge \neg F(s,\vec{s}) ,s\notin t(\vec{t})\vee F(s,\vec{t})}
}
{\vec{s}\neq \vec{t},(\exists x\in t(\vec{s}))\neg F(x,\vec{s}),s\notin t(\vec{t})\vee F(s,\vec{t})}
\end{align*}
If this proof tree is given by $\Phi_e(s)$, we set $\mathrm{Ext}(A,\vec{s},\vec{t})=\langle \code{b\forall},\code{(\forall x\in t(\vec{t}))F(x,\vec{t})},\code{\vec{s}\neq \vec{t},(\exists x\in t(\vec{s})\neg F(x,\vec{s}))},\mathrm{no}(\vec{s}\neq \vec{t},\neg A(\vec{s}),A(\vec{t})),0,e\rangle$. \par 
If $A$ is $\forall x\,F(x)$, $\mathrm{Ext}(A,\vec{s},\vec{t})$ is defined as follows. For every $s$, $\mathrm{Ext}(F(s,\vec{a}),\vec{s},\vec{t})$ derives $\vec{s}\neq \vec{t},\neg F(s,\vec{s}),F(s,\vec{t})$. Applying $\exists$ gives us $\vec{s}\neq \vec{t},\exists x\,\neg F(x,\vec{s}),F(s,\vec{t})$ for all $s$, and by $\forall$ we obtain $\vec{s}\neq \vec{t},\exists x\,\neg F(x,\vec{s}),\forall x\,F(x,\vec{t})$. \par 
The other cases are symmetrical. 
\end{proof}

\begin{sublemma}[Set induction]
KP proves the following. There is a partial recursive set function Ind such that, for any $\RS_\Omega(\bV)$-formula $F$ and $A\equiv \forall x((\forall y\in x)F(y)\to F(x))$, $\mathrm{Ind}(F)\in \RSderivation$, $\mathrm{End}(\mathrm{Ind}(F))=A\to \forall x\,F(x)$, $\mathrm{Length}(\mathrm{Ind}(F))\preceq \mathrm{no}(A\to \forall x\,F(x))\#\omega^{\rank(A)}$, $\mathrm{Rank}(\mathrm{Ind}(F))=0$. 
\end{sublemma}
\begin{proof}
We define a partial recursive set function $f$ in KP such that for any term $s$, $f(F,s)\in \RSderivation$, $\mathrm{End}(f(F,s))=\neg A,F(s)$, $\mathrm{Length}(f(F,s))\preceq \omega^{\rank(A)}\#\omega^{|s|+1}$, $\mathrm{Rank}(f(F,s))=0$; $f$ is defined by recursion on $|s|$. Given the index of $f(F,\cdot)$, it is then clear how to define Ind.  \par 
If $\sststile{0}{\omega^{\rank(A)}\#\omega^{|t|+1}}\neg A,F(t)$ has been proved by $f(F,t)$ for all $|t|<|s|$, then there is a partial recursive set function $g$ such that $g(F,t)$ proves $\sststile{0}{\omega^{\rank(A)}\#\omega^{|s|}+1}\neg A,t\in s\to F(t)$ for all terms $t$ (if $|t|<|s|$, this follows from $f(F,t)$; if $|t|\geq |s|$, then $t\notin s$ is an axiom). Now the derivation continues as follows:
\begin{align*}
\prftree
{\prftree
{\prftree
{\neg A,t\in s\to F(t)\text{ for all $t$}}
{\neg A,(\forall y\in s)F(y)}
\qquad
\neg F(s),F(s)
}
{\neg A,(\forall y\in s)F(y)\wedge \neg F(s),F(s)}
}
{\neg A,\exists x((\forall y\in x)F(y)\wedge \neg F(x) ),F(s)}
\end{align*}
where weakening is implied where appropriate; note that $\mathrm{no}(\neg F(s),F(s))<\omega^{\rank(A)}$, so the end sequent has length $\omega^{\rank(A)}\#\omega^{|s|}+4\prec \omega^{\rank(A)}\#\omega^{|s|+1}$; the end sequent is $\neg A,F(s)$ by contraction. 
\end{proof}

\begin{sublemma}[Pair]
KP proves the following. There is a partial recursive set function Pair such that, if $s,t$ are terms, then $\mathrm{Pair}(s,t)\in \RSderivation$, $\mathrm{End}(\mathrm{Pair}(s,t))=\exists z(s\in z\wedge t\in z)$, $\mathrm{Length}(\mathrm{Pair}(s,t))=1$, $\mathrm{Rank}(\mathrm{Pair}(s,t))=0$. 
\end{sublemma}
\begin{proof}
$s\in \{ s,t\}\wedge t\in \{ s,t\}$ is a $\Delta_0$-formula true in KP. 
\end{proof}

\begin{sublemma}[Union]
KP proves the following. There is a partial recursive set function Union such that, if $s$ is a term, then $\mathrm{Union}(s)\in \RSderivation$, $\mathrm{End}(\mathrm{Union}(s))=\exists z(\forall y\in s)(\forall x\in y)(x\in z)$, $\mathrm{Length}(\mathrm{Union}(s))=1$, $\mathrm{Rank}(\mathrm{Union}(s))=0$. 
\end{sublemma}
\begin{proof}
$(\forall y\in s)(\forall x\in y)x\in \bigcup s$ is a $\Delta_0$-formula true in KP. 
\end{proof}

\begin{sublemma}[Infinity]
KP proves the following. There is $w\in \RSderivation$ such that $\mathrm{End}(w)=\exists x( (\exists z\in x)z\in x\wedge (\forall y\in x)(\exists z\in x)y\in z )$ and $\mathrm{Length}(w)=1$, $\mathrm{Rank}(w)=0$; we also write $w$ as Inf.
\end{sublemma}
\begin{proof}
$(\exists z\in \omega)z\in \omega\wedge (\forall y\in \omega)(\exists z\in \omega)y\in z$ is a $\Delta_0$-formula true in KP. 
\end{proof}

\begin{sublemma}[$\Delta_0$-separation]
KP proves the following. There is a partial recursive set function Sep such that, given $A(a,b_1,\ldots,b_n)$ a $\Delta_0$-formula of KP with all free variables indicated, $s,t_1,\ldots,t_n$ terms, $\mathrm{Sep}(A,\langle s,\vec{t}\rangle )\in \RSderivation$, \\
$\mathrm{End}(\mathrm{Sep}(A, \langle s,\vec{t}\rangle))=\exists y((\forall x\in y)(x\in s)\wedge A(x,\vec{t})\wedge (\forall x\in s)(A(x,\vec{t})\to x\in y) )$, \\
$\mathrm{Length}(\mathrm{Sep}(A, \langle s,\vec{t}\rangle))=1$, and\\
$\mathrm{Rank}(\mathrm{Sep}(A, \langle s,\vec{t}\rangle))=0$. 
\end{sublemma}
\begin{proof}
$(\forall x\in \{ x\in s\mid A(x,\vec{t})\} )(x\in s\wedge A(x,\vec{t}))\wedge (\forall x\in s)(A(x,\vec{t})\to x\in \{ x\in s\mid  A(x,\vec{t})\} )$ is a $\Delta_0$-formula true in KP.
\end{proof}

\begin{sublemma}[$\Delta_0$-collection]
KP proves the following. There is a partial recursive set function Col such that, given $F(a,b,c_1,\ldots,c_n)$ a $\Delta_0$-formula of KP with all free variables indicated, $s,t_1,\ldots,t_n$ terms, \\
$\mathrm{Col}(F,\langle s,\vec{t}\rangle)\in \RSderivation$, \\
$\mathrm{End}(\mathrm{Col}(F,\langle s,\vec{t}\rangle))=(\forall x\in s)\exists y\,F(x,y,\vec{t})\to \exists z(\forall x\in s)(\exists y\in z)F(x,y,\vec{t})$, \\
$\mathrm{Length}(\mathrm{Col}(F,\langle s,\vec{t}\rangle))=\mathrm{no}((\forall x\in s)\exists y\,F(x,y,\vec{t})\to \exists z(\forall x\in s)(\exists y\in z)F(x,y,\vec{t}))$, \\
$\mathrm{Rank}(\mathrm{Col}(F,\langle s,\vec{t}\rangle))=0$. 
\end{sublemma}
\begin{proof}
By LEM we have 
\begin{align*}
\dststile{}{}\neg (\forall x\in s)\exists y\,F(x,y,\vec{t}),(\forall x\in s)\exists y\,F(x,y,\vec{t}).
\end{align*}
Applying $\Sigma$-Ref gives
\begin{align*}
\sststile{0}{\alpha+1}\neg(\forall x\in s)\exists y\,F(x,y,\vec{t}),\exists z(\forall x\in s)(\exists y\in z)F(x,y,\vec{t})
\end{align*}
where $\alpha=\mathrm{no}(\neg (\forall x\in s)\exists y\,F(x,y,\vec{t}),(\forall x\in s)\exists y\,F(x,y,\vec{t}))$. Two applications of $\vee$ gives
\begin{align*}
\sststile{0}{\alpha+3}(\forall x\in s)\exists y\,F(x,y,\vec{t})\to \exists z(\forall x\in s)(\exists y\in z)F(x,y,\vec{t});
\end{align*}
note also that $\alpha+3<\mathrm{no}((\forall x\in s)\exists y\,F(x,y,\vec{t})\to \exists z(\forall x\in s)(\exists y\in z)F(x,y,\vec{t}))$.
\end{proof}

Thus we have all the proofs for \fullref{embedding lemmas}.\par

Next we define a finitary proof system for KP.

\begin{named}[Definition of the finitary sequent calculus of $\KP$]
Let $a,b$ etc. denote free variables, $s,t$ etc. set terms, $\Gamma$ denote a finite set of formulas. KP has the following axioms:
\begin{itemize}
\item (Logical axioms). $\Gamma,A,\neg A$ for any formula $A$. 
\item (Extensionality). $\Gamma, a=b\wedge B(a)\to B(b)$ for any formula $B(a)$. 
\item (Set induction). $\Gamma,\forall x((\forall y\in x)F(y)\to F(x) )\to \forall x\,F(x)$ for any formula $F(a)$.
\item (Pair). $\Gamma,\exists z(a\in z\wedge b\in z)$. 
\item (Union). $\Gamma,\exists z(\forall y\in z)(\forall x\in y)x\in z$. 
\item (Infinity). $\Gamma,\exists x((\exists z\in x)z\in x\wedge (\forall y\in x)(\exists z\in x)y\in z)$. 
\item ($\Delta_0$-separation). $\Gamma,\exists y((\forall x\in y)(x\in a\wedge B(x))\wedge (\forall x\in a)(B(x)\to x\in y) )$ for any $\Delta_0$-formula $B$. 
\item ($\Delta_0$-collection). $\Gamma,(\forall x\in a)\exists y\,G(x,y)\to \exists z(\forall x\in a)(\exists y\in z)G(x,y)$ for any $\Delta_0$-formula $G$. 
\end{itemize}
The rules of inference are 
\begin{align*}
&(\wedge)\quad \prftree{ \Gamma,A\quad \Gamma,B}{\Gamma,A\wedge B}\\
&(\vee) \quad \prftree{ \Gamma,A}{\Gamma,A\vee B}\quad \prftree{ \Gamma,B}{\Gamma,A\vee B}\\
& (b\forall) \quad \prftree{\Gamma,a\in t\to F(a)}{\Gamma,(\forall x\in t)F(x)}  \\
& (b\exists) \quad \prftree{\Gamma,s\in t\wedge F(s)}{\Gamma,(\exists x\in t)F(x)}   \\
& (\forall)\quad \prftree{\Gamma,F(a)}{\Gamma,\forall x\,F(x)}\\
& (\exists )\quad \prftree{\Gamma,F(s)}{\Gamma,\exists x\,F(x)}   \\
& (\text{Cut}) \quad \prftree{\Gamma,A\quad \Gamma,\neg A}{\Gamma} 
\end{align*}
In inferring a sequent by $b\forall$ or $\forall$, the variable $a$ cannot appear in the lower sequents. \par 
We say that KP proves $\Gamma$ if $\Gamma$ can be derived from these axioms and inference rules. 
\end{named}

\begin{namedrmk*}[Notation]
Let $\phi_0$ denote the exponentiation $\alpha\mapsto \omega^\alpha$. 
\end{namedrmk*}

\begin{named}[Embedding theorem]
\label{embedding theorem}
KP proves the following. There is a partial recursive set function $P$ such that:
\begin{enumerate}
\item If $p$ is a code of an axiom of KP with end sequent $\Gamma(a_1,\ldots,a_n)$ where $a_1,\ldots,a_n$ are all the free variables that $\Gamma$ has, then for all terms $s_1,\ldots,s_n$, \\
$P(p,\langle s_1,\ldots,s_n\rangle)\in \RSderivation$, \\
$\mathrm{End}(P(p,\langle s_1,\ldots,s_n\rangle))=\Gamma(s_1,\ldots,s_n)$, \\
$\mathrm{Length}(P(p,\langle s_1,\ldots,s_n\rangle))\prec \Omega\cdot \omega^\omega$. 
\item If $p$ is a code of a KP-proof that is not an axiom and uses $k\geq 0$ instances of $b\forall/\forall$-inferences, with end sequent $\Gamma(a_1,\ldots,a_n)$ where $a_1,\ldots,a_n$ are all the free variables that $\Gamma$ has, then, there is some $m<\omega$ such that for all terms $s_1,\ldots,s_n$, \\
$P(p,\langle s_1,\ldots,s_n\rangle)\in \RSderivation$, \\
$\mathrm{End}(P(p,\langle s_1,\ldots,s_n\rangle))=\Gamma(s_1,\ldots,s_n)$, \\
$\mathrm{Length}(P(p,\langle s_1,\ldots,s_n\rangle))\prec \Omega\cdot \phi_0^{k+1}(\omega)$, \\
$\mathrm{Rank}(P(p,\langle s_1,\ldots,s_n\rangle))=\Omega+m$. 
\end{enumerate}
\end{named}
\begin{proof}
$P$ is defined by recursion on the length of $p$. \par 
If $p$ is an axiom, the value of $P$ is given by the previously defined functions LEM, Ext, Ind, Pair, Union, Inf, Sep, Col, along with appropriate weakening. Note that, if we write $\Gamma=A_1,\ldots,A_n$, by \ref{rank observation} (i), for some $m_1,\ldots,m_n$, we have $\mathrm{rank}(A_i)\preceq \omega\cdot \Omega +m_i$ for $i=1,\ldots,n$. Then $\mathrm{no}(\Gamma)=\omega^{\rank(A_1)}\#\ldots \#\omega^{\rank(A_n)}=(\Omega\cdot \omega^{m_1})\#\ldots \#(\Omega\cdot\omega^{m_n})=\Omega\cdot(\omega^{m_1}\#\ldots\#\omega^{m_n})\prec \Omega\cdot \omega^m$ for $m=\max(m_1,\ldots,m_n)+1$. \par 
Now assume $p$ is not an axiom. \par 
If the last inference of $p$ is $\wedge$, and we are in the situation 
\begin{align*}
\prftree{p_0\vdash\Gamma,A\quad p_1\vdash\Gamma,B}{p\vdash\Gamma,A\wedge B}
\end{align*}
then \\
$P(p,\langle s_1,\ldots,s_n\rangle)=\langle \code{\wedge},\code{A\wedge B},\code{\Gamma},a,r,P(p_0,\langle s_1,\ldots,s_n\rangle), P(p_1,\langle s_1,\ldots,s_n\rangle)\rangle$, where \\
$a=\max(\mathrm{Length}(P(p_0,\langle s_1,\ldots,s_n\rangle)), \mathrm{Length}(P(p_1,\langle s_1,\ldots,s_n\rangle)))+1$,\\
$r=\max(\Omega,\mathrm{Rank}(P(p_0,\langle s_1,\ldots,s_n\rangle)), \mathrm{Rank}(P(p_1,\langle s_1,\ldots,s_n\rangle)))$.  \par 
The cases $\vee$ and Cut are similarly dealt with. \par 
If the last inference of $p$ is $\exists$, and we are in the situation 
\begin{align*}
\prftree{p_0\vdash \Gamma,F(s(a_1,\ldots,a_n,b_1,\ldots,b_m))}{p\vdash\Gamma,\exists x\,F(x)}
\end{align*}
and $b_1,\ldots,b_m$ are free variables that do not appear in $\Gamma,\exists x\,F(x)$, then we can define \\
$P(p,\langle s_1,\ldots,s_n\rangle)=\langle \code{\exists},\code{\exists x\,F(x)},\code{\Gamma},a,r,P(p_0,\langle s_1,\ldots,s_n,c_\emptyset,\ldots,c_\emptyset\rangle)\rangle$, where \\
$a=\mathrm{Length}(P(p_0,\langle s_1,\ldots,s_n,c_\emptyset,\ldots,c_\emptyset\rangle))+1$,\\
$r=\max(\Omega,\mathrm{Rank}(P(p_0,\langle s_1,\ldots,s_n,c_\emptyset,\ldots,c_\emptyset\rangle)))$. \par 
The case $b\exists$ is similarly dealt with. \par 
If the last inference of $p$ is $\forall$, and we are in the situation
\begin{align*}
\prftree{p_0\vdash\Gamma,F(a)}{p\vdash\Gamma,\forall x\,F(x)}
\end{align*}
and so $p_0$ uses $k-1$ instances of $b\forall/\forall$-inferences, if $p$ uses $k$ of them. By induction, we have some $m<\omega$ such that for all terms $s_1,\ldots,s_n,t$,\\
$\mathrm{End}(P(p_0,\langle s_1,\ldots,s_n,t\rangle))=\Gamma(s_1,\ldots,s_n),F(s_1,\ldots,s_n,t)$, \\
$\mathrm{Length}(P(p_0,\langle s_1,\ldots,s_n,t\rangle))\prec \Omega\cdot\phi_0^k(\omega)$,\\
$\mathrm{Rank}(P(p_0,\langle s_1,\ldots,s_n,t\rangle))=\Omega+m$. \\
Therefore, if $e$ is an index for $\lambda t.P(p_0,\langle s_1,\ldots,s_n,t\rangle)$, we may define \\
$P(p,\langle s_1,\ldots,s_n\rangle)=\langle \code{\forall},\code{\forall x\,F(x)},\code{\Gamma},\Omega\cdot\phi_0^k(\omega),\mathrm{Rank}(P(p_0,\langle s_1,\ldots,s_n,c_\emptyset\rangle)), e\rangle$.\footnote{The uniformity of the ranks of the premises allows us to use $\mathrm{Rank}(P(p_0,\langle s_1,\ldots,s_n,c_\emptyset\rangle))$ as the rank of the derivation $P(p,\langle s_1,\ldots,s_n\rangle)$.}\par 
If the last inference of $p$ is $b\forall$, and we are in the situation
\begin{align*}
\prftree{p_0\vdash \Gamma(\vec{s}),a\in t(\vec{s})\to F(\vec{s}, a)}{p\vdash \Gamma(\vec{s}),(\forall x\in t(\vec{s}))F(\vec{s},x)}
\end{align*}
and $p_0$ uses $k-1$ instances of $b\forall/\forall$-inferences if $p$ uses $k$ of them. By induction, there is some $m<\omega$ such that for all terms $\vec{s},r$, \\
$\mathrm{End}(P(p_0,\langle \vec{s},r\rangle))=\Gamma(\vec{s}),r\in t(\vec{s})\to F(\vec{s},r)$,\\
$\mathrm{Length}(P(p_0,\langle \vec{s},r\rangle))\prec \Omega\cdot\phi_0^k(\omega)$,\\
$\mathrm{Rank}(P(p_0,\langle \vec{s},r\rangle))=\Omega+m$.\\ Therefore, if $e(\vec{s})$ is an index (function) for $\lambda r.P(p_0,\langle \vec{s},r\rangle)$, we may define \\
$P(p,\vec{s})=\langle \code{b\forall},\code{(\forall x\in t)F(x)},\code{\Gamma},\Omega\cdot\phi_0^k(\omega),\mathrm{Rank}(P(p_0,\langle \vec{s},c_\emptyset\rangle)),e(\vec{s})\rangle$.  
\end{proof}

\begin{cor}
For any finite set $\Gamma$ of formulas, $\KP+\TI(\epsilon_{\Omega+1})$ proves that, for all sets $x_1,\ldots,x_n$, if $\Gamma(c_{x_1},\ldots,c_{x_n})$ is KP-provable, then there is $w\in \RSderivation$ with end sequent $\Gamma(c_{x_1},\ldots,c_{x_n})$, $\mathrm{Length}(w)\prec \epsilon_{\Omega+1}$, $\mathrm{Rank}(w)\preceq \Omega+1$. 
\end{cor}
\begin{proof}
Given a KP-proof $p$ of $\Gamma(c_{x_1},\ldots,c_{x_n})$, we have $P(p,\langle\rangle)\in \RSderivation$, $\mathrm{End}(P(p,\langle\rangle))=\Gamma(c_{x_1},\ldots,c_{x_n})$, $\mathrm{Length}(P(p,\langle\rangle))\prec \Omega\cdot\epsilon_0$, $\mathrm{Rank}(P(p,\langle\rangle))=\Omega+m$ for some $m$ which we can read off from $P(p,\langle\rangle)$. Thus $w:=\mathrm{CutElim}^{m-1}(P(p,\langle\rangle))\in \RSderivation$ has the same end sequent as $P(p,\langle\rangle)$, the length of $w$ is $\preceq \phi_0^{m-1}(\Omega\cdot\epsilon_0)\prec \epsilon_{\Omega+1}$, and $\mathrm{Rank}(w)=\Omega+1$. 
\end{proof}

\begin{lem}
\label{induction along proofs with cuts eliminated}
$\KP+\TI(\epsilon_{\Omega+1})$ proves that, for any finite set $\Gamma$ of $\Pi_n$-formulas ($n\geq 2$), for all sets $x_1,\ldots,x_n$, if $w\in \RSderivation$ has end sequent $\Gamma(c_{x_1},\ldots,c_{x_n})$, $\mathrm{Length}(w)\prec \epsilon_{\Omega+1}$, $\mathrm{Rank}(w)\preceq \Omega+1$, then $\mathrm{True}_{\Pi_n}(\code{\bigvee\Gamma(x_1,\ldots,x_n)})$. 
\end{lem}
\begin{proof}
By induction on $a$ along $\epsilon_{\Omega+1}$. \par 
If $w=\langle \code{\wedge},\code{A\wedge B},\code{\Gamma},a,r,u,v\rangle$, by induction, the disjunctions of $\mathrm{End}(u)=\Gamma,A$ and $\mathrm{End}(v)=\Gamma,B$ are true, so the disjunction of $\Gamma,A\wedge B$ is true. \par 
If $w=\langle \code{\forall},\code{\forall x\,F(x)},\code{\Gamma},a,r,e\rangle$, by induction the disjunctions of $\mathrm{End}(\Phi_e(s))=\Gamma,F(s)$ are true for all $s$, so the disjunction of $\Gamma,\forall x\,F(x)$ is true. \par 
The other cases of logical rules are similar.\par  
If the last inference is a cut, then the cut formula must be $\Sigma_1$/$\Pi_1$ and so the equivalence $\True_{\Pi_n}(\code{\bigvee\Gamma})\leftrightarrow \bigvee\Gamma$ is not violated. If the last inference is $\Sigma$-Ref, the induction is completed by the fact that KP proves $A\leftrightarrow \exists z\,A^z$ for any $\Sigma$-formula $A$.  
\end{proof}

This concludes the proof of our theorem:

\begin{thm}
$\KP+\TI(\epsilon_{\Omega+1})\vdash \RFN(\KP)$. \qedbyhand 
\end{thm}

Since the proof of \fullref{induction along proofs with cuts eliminated} simply hinges on the fact that we may eliminate cuts above certain complexity that is only determined by the complexity of axioms in $\RS_\Omega(\bV)$, we may deduce that, if $T$ stands for $\KP+\Gamma\text{-separation}+\Gamma\text{-collection}$ where $\Gamma$ is any given syntactic complexity,
\begin{align*}
\RFN(T)\equiv  \TI(\epsilon_{\Omega+1})\quad\text{over }T
\end{align*}
as well; we may add the axiom schemata $\Gamma$-separation and $\Gamma$-collection into the system $\RS_\Omega(\bV)$, forgo the comprehension terms $\{ x\in s\mid A(x,\vec{t})\}$, and the cut elimination (now above the complexity of $\Gamma$-separation and $\Gamma$-collection) and the embedding theorem go through as usual. \par 

This argument does not apply to ZF; indeed, we can show that ZF proves $\TI(\epsilon_{\Omega+1})$. Suppose it doesn't, then there is a formula $F$ and a model $M$ of ZF such that $\TI(\epsilon_{\Omega+1},F)$ is false in $M$. By the reflection theorem of ZF, there is a set model $N$ in $M$ such that $\TI(\epsilon_{\Omega+1},F)$ is absolute for $N$, and so $\neg\TI(\epsilon_{\Omega+1},F)^N$. But as $(\epsilon_{\Omega+1})^N$ is an ordinal in $M$, this contradicts foundation in $M$.


\printbibliography[ heading=bibintoc, title={References}]

@book{Barwise1975,
  title={Admissible Sets and Structures},
  author={Barwise, J.},
  isbn={9781107168336},
  series={Perspectives in Logic},
  url={https://books.google.co.uk/books?id=3aYoDgAAQBAJ},
  year={1975},
  publisher={Cambridge University Press}
}

@article{Rathjen1992,
 ISSN = {00224812},
 URL = {http://www.jstor.org/stable/2275441},
 author = {Michael Rathjen},
 journal = {The Journal of Symbolic Logic},
 number = {3},
 pages = {954--969},
 publisher = {Association for Symbolic Logic},
 title = {A Proof-Theoretic Characterization of the Primitive Recursive Set Functions},
 volume = {57},
 year = {1992}
}

@incollection{Schwichtenberg1977,
title = {Proof Theory: Some Applications of Cut-Elimination},
editor = {Jon Barwise},
series = {Studies in Logic and the Foundations of Mathematics},
publisher = {Elsevier},
volume = {90},
pages = {867-895},
year = {1977},
booktitle = {HANDBOOK OF MATHEMATICAL LOGIC},
issn = {0049-237X},
doi = {https://doi.org/10.1016/S0049-237X(08)71124-8},
url = {https://www.sciencedirect.com/science/article/pii/S0049237X08711248},
author = {Helmut Schwichtenberg},
}

@article{KreiselLevy1968,
author = {Kreisel, G. and Lévy, A.},
title = {Reflection Principles and their Use for Establishing the Complexity of Axiomatic Systems},
journal = {Mathematical Logic Quarterly},
volume = {14},
number = {7-12},
pages = {97-142},
doi = {https://doi.org/10.1002/malq.19680140702},
url = {https://onlinelibrary.wiley.com/doi/abs/10.1002/malq.19680140702},
eprint = {https://onlinelibrary.wiley.com/doi/pdf/10.1002/malq.19680140702},
year = {1968}
}

@article{AczelRathjen2001,
author = {Peter Aczel and Michael Rathjen},
title = {Notes on constructive set theory},
year = {2001},
publisher = {Institut Mittag-Leffler},
url = {http://www.ml.kva.se/preprints/archive2000-2001.php}
}

@incollection{Normann1978,
title = {Set Recursion},
editor = {J.E. Fenstad and R.O. Gandy and G.E. Sacks},
series = {Studies in Logic and the Foundations of Mathematics},
publisher = {Elsevier},
volume = {94},
pages = {303-320},
year = {1978},
booktitle = {Generalized Recursion Theory II},
issn = {0049-237X},
doi = {https://doi.org/10.1016/S0049-237X(08)70938-8},
url = {https://www.sciencedirect.com/science/article/pii/S0049237X08709388},
author = {Dag Normann}
}

@book{Sacks2017, place={Cambridge}, series={Perspectives in Logic}, title={Higher Recursion Theory}, DOI={10.1017/9781316717301}, publisher={Cambridge University Press}, author={Sacks, Gerald E.}, year={2017}, collection={Perspectives in Logic}}

@article{Rathjen2012,
title = {From the weak to the strong existence property},
journal = {Annals of Pure and Applied Logic},
volume = {163},
number = {10},
pages = {1400-1418},
year = {2012},
note = {Set Theory, Classical and Constructive – Invited papers from the meeting in Amsterdam, May 6–7, 2010},
issn = {0168-0072},
doi = {https://doi.org/10.1016/j.apal.2012.01.012},
url = {https://www.sciencedirect.com/science/article/pii/S0168007212000310},
author = {Michael Rathjen}
}

@article{Parsons1970,
  title={On a Number Theoretic Choice Schema and its Relation to Induction},
  author={Charles D. Parsons},
  journal={Studies in logic and the foundations of mathematics},
  year={1970},
  volume={60},
  pages={459-473}
}

\end{document}